# Limit laws for the energy of a charged polymer

Xia Chen[1]

*Department of Mathematics, University of Tennessee, Knoxville, TN 37996-1300, USA. E-mail: xchen@math.utk.edu*

**Abstract.** In this paper we obtain the central limit theorems, moderate deviations and the laws of the iterated logarithm for the energy

$$H_n = \sum_{1\le j<k\le n} \omega_j\omega_k 1_{\{S_j=S_k\}}$$

of the polymer $\{S_1,\dots,S_n\}$ equipped with random electrical charges $\{\omega_1,\dots,\omega_n\}$. Our approach is based on comparison of the moments between $H_n$ and the self-intersection local time

$$Q_n = \sum_{1\le j<k\le n} 1_{\{S_j=S_k\}}$$

run by the $d$-dimensional random walk $\{S_k\}$. As partially needed for our main objective and partially motivated by their independent interest, the central limit theorems and exponential integrability for $Q_n$ are also investigated in the case $d\ge 3$.

**Résumé.** Cet article est consacré à l'étude du théorème central limite, des déviations modérées et des lois du logarithme itéré pour l'énergie

$$H_n = \sum_{1\le j<k\le n} \omega_j\omega_k 1_{\{S_j=S_k\}}$$

du polymère $\{S_1,\dots,S_n\}$ doté de charges électriques $\{\omega_1,\dots,\omega_n\}$. Notre approche se base sur la comparaison des moments de $H_n$ et du temps local de recoupements

$$Q_n = \sum_{1\le j<k\le n} 1_{\{S_j=S_k\}}$$

de la marche aléatoire $d$-dimensionelle $\{S_k\}$. L'étude du théorème central limite et de l'intégrabilité exponentielle de $Q_n$ (dans le cas $d\ge 3$) est également menée, tant pour comme outil pour notre principal objectif que pour son intérêt intrinsèque.

[1]Research partially supported by NSF Grant DMS-0704024.

## 1. Introduction

In the physics literature, the geometric shape of certain polymers is often described as an interpolation line segment with the vertices given as the $n$-step lattice (simple) random walk

$$\{S_1, S_2, \ldots, S_n\}.$$

By placing independent, identically distributed electric charges $\omega_k = \pm 1$ to each vertex of the polymer, Kantor and Kardar [16] consider a model of polymers with random electrical charges associated with the Hamiltonian

$$H_n = \sum_{1 \le j < k \le n} \omega_j \omega_k 1_{\{S_j = S_k\}}. \tag{1.1}$$

In the physics literature, $H_n$ is called the energy of the polymer. To understand the physics intuition of $H_n$, we assign an electrical charge $\omega_k$ to the random site $S_k$ for all $k = 1, 2, \ldots$. Assume that when two charges meet, the pair with opposite signs gives negative contribution while the pair with the same sign gives positive contribution. Thus, $H_n$ represents the total electrical interaction charge of the polymer $\{S_1, S_2, \ldots, S_n\}$.

We point out some other works by physicists in this direction. In [10], the charges are i.i.d. Gaussian variables. In [11], the charges take 0–1 values. We also refer the reader to [4, 18] for the continuous versions of the polymer with random charges. Finally, we mention the survey paper by van der Hofstad and König [12] for a long list of mathematical models connected to polymers.

As for other connections, we cite the comment by Martínez and Petritis [18]: "It is argued that a protein molecule is very much like a random walk with random charges attached at the vertices of the walk; these charges are interacting through local interactions mimicking Lennard–Jones or hydrogen-bond potentials".

We study the asymptotic behaviors of $H_n$ given in (1.1). In the rest of the paper, $\{S_n\}_{n\ge1}$ is a symmetric random walk on $\mathbb{Z}^d$ with covariance matrix $\Gamma$ (or variance $\sigma^2$ as $d=1$). We assume that the smallest group that supports $\{S_n\}_{n\ge1}$ is $\mathbb{Z}^d$. Throughout, $\{\omega_k\}_{k\ge1}$ is an i.i.d. sequence of symmetric random variable with

$$\mathbb{E}\omega_1^2 = 1 \quad \text{and} \quad \mathbb{E}e^{\lambda_0 \omega_1^2} < \infty \quad \text{for some } \lambda_0 > 0. \tag{1.2}$$

Our first result is on the central limit theorems.

**Theorem 1.1.** *As $d=1$,*

$$\frac{1}{n^{3/4}} H_n \xrightarrow{d} (2\sigma)^{-1/2} \left( \int_{-\infty}^{\infty} L^2(1,x)\, dx \right)^{1/2} U, \tag{1.3}$$

*where $U$ is a random variable with standard normal distribution, $L(t,x)$ is the local time of the 1-dimensional Brownian motion $W(t)$ such that $U$ and $W(t)$ are independent.*

*As $d=2$,*

$$\frac{1}{\sqrt{n \log n}} H_n \xrightarrow{d} \frac{1}{\sqrt{2\pi \sqrt[4]{\det \Gamma}}} U. \tag{1.4}$$

*As $d \ge 3$,*

$$\frac{1}{\sqrt{n}} H_n \xrightarrow{d} \sqrt{\gamma} U, \tag{1.5}$$

*where*

$$\gamma = \sum_{k=1}^{\infty} \mathbb{P}\{S_k = 0\}. \tag{1.6}$$

Here is our explanation on the dimensional dependence appearing in Theorem 1.1. The higher the dimension is, the less likely the random walk is to have long-range interaction (self-intersection). In the multi-dimensional case ($d \geq 2$), therefore, $H_n$ is a sum of random variables with weak dependence and yields a Gaussian limit when properly normalized. It should be pointed that the low level of long-range interaction is vital for the chaos

$$\sum_{1\leq j<k\leq n} a_{j,k}\omega_j\omega_k$$

to have a Gaussian limit when properly normalized. A simple example is when $a_{j,k} \equiv 1$. In this case

$$\sum_{1\leq j<k\leq n} \omega_j\omega_k = \frac{1}{2}\left\{\left[\sum_{j=1}^n \omega_j\right]^2 - \sum_{j=1}^n \omega_j^2\right\}.$$

By the classic law of large numbers and classic central limit theorem,

$$\frac{1}{n}\sum_{1\leq j<k\leq n} \omega_j\omega_k \xrightarrow{d} \frac{1}{2}(U^2-1)$$

which sharply contrasts to the statements in Theorem 1.1.

Our next theorem describes the moderate deviation behaviors of $H_n$.

**Theorem 1.2.** *As* $d = 1$,

$$\lim_{n\to\infty} \frac{1}{b_n}\log \mathbb{P}\{\pm H_n \geq \lambda(nb_n)^{3/4}\} = -\frac{1}{2}\sigma^{2/3}(3\lambda)^{4/3}, \quad \lambda > 0, \tag{1.7}$$

*for any positive sequence* $\{b_n\}$ *satisfying*

$$b_n \to \infty \quad \text{and} \quad b_n = \mathrm{o}(\sqrt[7]{n}), \quad n \to \infty. \tag{1.8}$$

*As* $d = 2$,

$$\lim_{n\to\infty} \frac{1}{b_n}\log \mathbb{P}\{\pm H_n \geq \lambda\sqrt{n(\log n)b_n}\} = -\pi\sqrt{\det(\Gamma)}\lambda^2, \quad \lambda > 0, \tag{1.9}$$

*for any positive sequence* $\{b_n\}$ *satisfying*

$$b_n \to \infty \quad \text{and} \quad b_n = \mathrm{o}(\log n), \quad n \to \infty. \tag{1.10}$$

*As* $d \geq 3$,

$$\lim_{n\to\infty} \frac{1}{b_n}\log \mathbb{P}\{\pm H_n \geq \lambda\sqrt{nb_n}\} = -\frac{\lambda^2}{2\gamma}, \quad \lambda > 0, \tag{1.11}$$

*for any positive sequence* $\{b_n\}$ *satisfying*

$$b_n \to \infty \quad \text{and} \quad b_n = \mathrm{o}\left(\left(\frac{n}{\log n}\right)^{1/4}\right), \quad n \to \infty. \tag{1.12}$$

Our moderate deviations applied to the law of the iterated logarithm:

**Theorem 1.3.** *As* $d = 1$,

$$\limsup_{n\to\infty} \frac{\pm H_n}{(n\log\log n)^{3/4}} = \frac{2^{3/4}}{3}\sigma^{-1/2} \quad a.s. \tag{1.13}$$

*As* $d=2$,

$$\limsup_{n\to\infty} \frac{\pm H_n}{\sqrt{n\log n\log\log n}} = \frac{1}{\sqrt{\pi}\det(\Gamma)^{1/4}} \quad a.s. \tag{1.14}$$

*As* $d\ge 3$,

$$\limsup_{n\to\infty} \frac{\pm H_n}{\sqrt{n\log\log n}} = \sqrt{2\gamma} \quad a.s. \tag{1.15}$$

We compare the results and treatments between the present paper and some recent works on self-intersection local times such as [3, 6]. On the one hand, we shall see that the asymptotic behaviors of $H_n$ described in our main theorems are closely related to those of the self-intersection local time

$$Q_n = \sum_{1\le j<k\le n} 1_{\{S_j=S_k\}}. \tag{1.16}$$

Indeed, our approach is based on the moment comparisons between $H_n$ and $Q_n$ (see Proposition 2.1). In particular, the difference in limiting distribution between the case $d=1$ and the case $d\ge 2$ in Theorem 1.1 is caused by the fact that in the case $d=1$, $Q_n$ converges (in distribution) to the Brownian self-intersection local times when properly normalized (see [8]), while as $d\ge 2$, $Q_n$ is asymptotically close to its expectation (see [3] for $d=2$ and the Section 5 for $d\ge 3$).

On the other hand, the fact that $Q_n$ is close to the quadratic form

$$\sum_{x\in\mathbb{Z}^d} l^2(n,x)$$

of the local time $l(n,x)$ plays a crucial role in the study of the self-intersection local time $Q_n$ (see e.g., [3, 8]). It allows, for example, some technologies developed along the line of probability in Banach space. Unfortunately, this idea does not work in our setting. Indeed, the fact (in view of Theorem 1.1) that the second term in the decomposition

$$\sum_{x\in\mathbb{Z}^d}\left[\sum_{j=1}^n \omega_j 1_{\{S_j=x\}}\right]^2 = 2H_n + \sum_{j=1}^n \omega_j^2 \tag{1.17}$$

is the dominating term shows that $H_n$ is not even in the same asymptotic order as the quadratic form on the left-hand side.

The key estimations are carried out in Proposition 2.1. Our approach relies on the following crucial observation. By (1.17) we have

$$H_n = \frac{1}{2}\sum_{x\in\mathbb{Z}^d}\left\{\left[\sum_{j=1}^n \omega_j 1_{\{S_j=x\}}\right]^2 - \sum_{j=1}^n \omega_j 1_{\{S_j=x\}}\right\}. \tag{1.18}$$

Conditioned on the random walk $\{S_k\}$, the random variables

$$\left[\sum_{j=1}^n \omega_j 1_{\{S_j=x\}}\right]^2 - \sum_{j=1}^n \omega_j 1_{\{S_j=x\}}, \quad x\in\mathbb{Z}^d,$$

form an independent family and, for each fix $x\in\mathbb{Z}^d$,

$$\left[\sum_{j=1}^n \omega_j 1_{\{S_j=x\}}\right]^2 - \sum_{j=1}^n \omega_j 1_{\{S_j=x\}} \overset{d}{=} \left[\sum_{j=1}^{l(n,x)} \omega_j\right]^2 - \sum_{j=1}^{l(n,x)} \omega_j.$$

By a classic estimate for independent sums, and by some combinatorial computation, a conditional moment estimate given in Proposition 2.1 links $H_n$ with $Q_n$. Another fact repeatedly used in this work is that the self-intersection occurring at the frequently visited sites does not make a significant contribution to the quantities $H_n$ and $Q_n$. Consequently, the pairs $H_n$ and $\widetilde{H}_n$ (defined in (2.3)); $Q_n$ and $\widetilde{Q}_n$ (defined in (2.4) below) are exchangeable in our setting.

Beyond mathematical technicality, the creation of the present paper is based on our belief that $H_n$ resembles, in the limiting behaviors described in our main theorems, the random quantity

$$H_n' = \sum_{1\le j<k\le n} 1_{\{S_j=S_k\}} U_{j,k},$$

where $U_{j,k}$ are i.i.d. standard normal random variables independent of $\{S_k\}$. Notice that $H_n'$ is conditionally normal with conditional variance $Q_n$. Our observation explains why and how the limiting behaviors of $H_n$ depend on its conditional variance $Q_n$. It should be pointed out, however, that the replacement of $\omega_j\omega_k$ by $U_{j,k}$ is highly non-trivial and should not be taken for granted, in view of the example given next to Theorem 1.1.

In Section 2, we establish a comparison (Proposition 2.1) between the moments of $H_n$ and $Q_n$, and then apply it to prove Theorem 1.1. Our approach relies on combinatorial and conditioning methods. In Section 3, Proposition 2.1 is further applied to prove Theorem 1.2 through a Laplacian argument. In Section 4, the laws of the iterated logarithm given in Theorem 1.3 are proved as a consequence of our moderate deviations. The non-trivial part of this section is a maximal inequality (Lemma 4.1) of Lévy type. In Section 5, we investigate the weak laws and exponential integrabilities for the renormalized self-intersection local time $Q_n - \mathbb{E}Q_n$ in the high dimensions ($d\ge 3$). The central limit theorem given in Theorem 5.1 and the exponential integrability given in Theorem 5.2 provide sharp bounds on $Q_n - \mathbb{E}Q_n$, which constitute the replacement of $Q_n$ by $\mathbb{E}Q_n$ carried out in our argument for Theorem 1.1 and for Theorem 1.2 (the estimate of $Q_n - \mathbb{E}Q_n$ needed in the case $d=2$ was established in [3, 20]). In addition, the results given in Section 5 are of independent interest as a part of the study of the self-intersection local times in high dimensions and are partially motivated by some recent works of Asselah and Castell [1] and Asselah [2].

## 2. Moment comparison and laws of weak convergence

We begin with the following classic lemma.

**Lemma 2.1.** *Assume (1.2). Then*

$$\mathbb{E}\left\{\left(\sum_{j=1}^n \omega_j\right)^2 - \sum_{j=1}^n \omega_j^2\right\}^2 = 2n(n-1). \tag{2.1}$$

*More generally, there is a constant $C>0$ such that for any integers $n\ge 1$ and $m\ge 2$,*

$$\mathbb{E}\left|\left(\sum_{j=1}^n \omega_j\right)^2 - \sum_{j=1}^n \omega_j^2\right|^m \le m!(Cn(n-1))^{m/2}. \tag{2.2}$$

**Proof.** The first part follows from the following straightforward computation:

$$\mathbb{E}\left\{\left(\sum_{j=1}^n \omega_j\right)^2 - \sum_{j=1}^n \omega_j^2\right\}^2 = 4\mathbb{E}\left\{\sum_{1\le j<k\le n} \omega_j\omega_k\right\}^2 = 4\sum_{1\le j<k\le n} \mathbb{E}(\omega_j^2\omega_k^2) = 2n(n-1).$$

For the second part, we only need to show

$$\mathbb{E}\left|\left(\sum_{j=1}^n \omega_j\right)^2 - \sum_{j=1}^n \omega_j^2\right|^m \le m! C^{m/2} n^m.$$

By the inequality

$$\left(\mathbb{E}\left|\left(\sum_{j=1}^n\omega_j\right)^2-\sum_{j=1}^n\omega_j^2\right|^m\right)^{1/m}\le\left(\mathbb{E}\left(\sum_{j=1}^n\omega_j\right)^{2m}\right)^{1/m}+\left(\mathbb{E}\left(\sum_{j=1}^n\omega_j^2\right)^m\right)^{1/m}$$

all we need is that

$$\mathbb{E}\left(\sum_{j=1}^n\omega_j\right)^{2m}\le C^{m/2}m!n^m$$

and that

$$\mathbb{E}\left|\sum_{j=1}^n\omega_j^2\right|^m\le C^{m/2}m!n^m.$$

Due to similarity we only prove the first inequality. Notice that by symmetry

$$\mathbb{E}\left(\sum_{j=1}^n\omega_j\right)^{2m}=\sum_{\substack{k_1+\cdots+k_n=m\\ k_1,\dots,k_n\ge 0}}\frac{(2m)!}{(2k_1)!\cdots(2k_n)!}\mathbb{E}\omega^{2k_1}\cdots\mathbb{E}\omega^{2k_n}.$$

By the integrability given in (1.2) there is a constant $c_1>0$ such that

$$\mathbb{E}\omega^{2k}\le k!c_1^k,\quad k=0,1,2,\dots.$$

Notice also the very rough estimate

$$c_2^k(k!)^2\le(2k)!\le c_3^k(k!)^2,\quad k=0,1,2,\dots.$$

So we have

$$\mathbb{E}\left(\sum_{j=1}^n\omega_j\right)^{2m}\le C^{m/2}m!\sum_{\substack{k_1+\cdots+k_n=m\\ k_1,\dots,k_n\ge 0}}\frac{m!}{k_1!\cdots k_n!}=C^{m/2}m!n^m.$$

□

Let $K_n$ be a positive sequence which may vary in different settings and will later be specified in each specific setting. Recall that $Q_n$ is given in (1.16) and define the local time

$$l(n,x)=\sum_{k=1}^n 1_{\{S_k=x\}},\quad x\in\mathbb{Z}^d, n=1,2,\dots.$$

The asymptotic behaviors of the local times of the random walks have been studied extensively. We cite the book by Révész [19] for an overview.

The following two random quantities play important roles in this paper:

$$\widetilde{H}_n=H_n 1_{\{\sup_{x\in\mathbb{Z}^d}l(n,x)\le K_n\}},\tag{2.3}$$

$$\widetilde{Q}_n=Q_n 1_{\{\sup_{x\in\mathbb{Z}^d}l(n,x)\le K_n\}}.\tag{2.4}$$

In addition, we introduce the deterministic quantity

$$A_m(n)=\frac{1}{2^m}\sum_{(y_1,\dots,y_m)\in B_m}\mathbb{E}\left(1_{\{\sup_{x\in\mathbb{Z}^d}l(n,x)\le K_n\}}\prod_{k=1}^m l(n,y_k)(l(n,y_k)-1)\right),$$

where $m,n=1,2,\dots$ and

$$B_m=\{(y_1,\dots,y_m)\in(\mathbb{Z}^d)^m; y_1,\dots,y_m \text{ are distinct}\}. \tag{2.5}$$

An easy observation gives that

$$\begin{aligned} A_m(n) &\le \frac{1}{2^m}\sum_{y_1,\dots,y_m\in\mathbb{Z}^d}\mathbb{E}\left(1_{\{\sup_{x\in\mathbb{Z}^d} l(n,x)\le K_n\}}\prod_{k=1}^m l(n,y_k)(l(n,y_k)-1)\right)\\ &=\mathbb{E}\widetilde{Q}_n^m. \end{aligned} \tag{2.6}$$

Some more substantial comparisons are given in the following.

**Proposition 2.1.** *There is a constant $C>0$ independent of $n$, $m$ and the choice of $K_n$, such that*

$$\mathbb{E}\widetilde{H}_n^m\le\frac{m!}{2^m}\sum_{l=1}^{[2^{-1}m]}\frac{1}{l!}K_n^{m-2l}2^lC^{(m-2l)/2}\binom{m-l-1}{m-2l}\mathbb{E}\widetilde{Q}_n^l. \tag{2.7}$$

*On the other hand, for any integers $m,n\ge1$,*

$$\mathbb{E}\widetilde{H}_n^{2m}\ge\frac{(2m)!}{2^m m!}A_m(n), \tag{2.8}$$

$$\mathbb{E}\widetilde{Q}_n^m\le\sum_{l=1}^m\binom{m}{l}\left(\frac{lK_n^2}{2}\right)^{m-l}A_l(n). \tag{2.9}$$

**Proof.** Notice that

$$H_n=\frac{1}{2}\sum_{x\in\mathbb{Z}^d}\left\{\left(\sum_{j=1}^n\omega_j1_{\{S_j=x\}}\right)^2-\sum_{j=1}^n\omega_j^21_{\{S_j=x\}}\right\}=\frac{1}{2}\sum_{x\in\mathbb{Z}^d}\Lambda_n(x)\qquad\text{(say)}. \tag{2.10}$$

Hence,

$$\mathbb{E}\widetilde{H}_n^m=2^{-m}\sum_{x_1,\dots,x_m\in\mathbb{Z}^d}\mathbb{E}\left(1_{\{\sup_{x\in\mathbb{Z}^d} l(n,x)\le K_n\}}\prod_{k=1}^m\Lambda_n(x_k)\right).$$

For each $1\le l\le m$, let

$$A_l=\{(x_1,\dots,x_m)\in(\mathbb{Z}^d)^m;\#\{x_1,\dots,x_m\}=l\}.$$

Then,

$$\mathbb{E}\widetilde{H}_n^m=2^{-m}\sum_{l=1}^m\sum_{(x_1,\dots,x_m)\in A_l}\mathbb{E}\left(1_{\{\sup_{x\in\mathbb{Z}^d} l(n,x)\le K_n\}}\prod_{k=1}^m\Lambda_n(x_k)\right). \tag{2.11}$$

Write

$$\mathcal{C}_l=\{F\subset\mathbb{Z}^d;\#(F)=l\}$$

and for any $\{y_1,\dots,y_l\}\in\mathcal{C}_l$, set

$$A_l(y_1,\dots,y_l)=\{(x_1,\dots,x_m)\in(\mathbb{Z}^d)^m;\{x_1,\dots,x_m\}=\{y_1,\dots,y_l\}\}.$$

Notice that

$$1_{A_l}(x_1,\dots,x_m) = \sum_{\{y_1,\dots,y_l\}\in\mathcal{C}_l} 1_{A_l(y_1,\dots,y_l)}(x_1,\dots,x_m).$$

Thus

$$\sum_{(x_1,\dots,x_m)\in A_l} \mathbb{E}\Bigg(1_{\{\sup_{x\in\mathbb{Z}^d} l(n,x)\le K_n\}} \prod_{k=1}^m \varLambda_n(x_k)\Bigg)$$
$$= \sum_{\{y_1,\dots,y_l\}\in\mathcal{C}_l} \sum_{(x_1,\dots,x_m)\in A_l(y_1,\dots,y_l)} \mathbb{E}\Bigg(1_{\{\sup_{x\in\mathbb{Z}^d} l(n,x)\le K_n\}} \prod_{k=1}^m \varLambda_n(x_k)\Bigg).$$

For any $(x_1,\dots,x_m)\in A_l(y_1,\dots,y_l)$, let $i_k$ be the number of $x_1,\dots,x_m$ which are equal to $y_k$, where $k=1,\dots,l$. Then

$$\mathbb{E}\Bigg(1_{\{\sup_{x\in\mathbb{Z}^d} l(n,x)\le K_n\}} \prod_{k=1}^m \varLambda_n(x_k)\Bigg) = \mathbb{E}\Bigg\{1_{\{\sup_{x\in\mathbb{Z}^d} l(n,x)\le K_n\}} \prod_{k=1}^l \varLambda_n(y_k)^{i_k}\Bigg\}.$$

Consequently,

$$\sum_{(x_1,\dots,x_m)\in A_l(y_1,\dots,y_l)} \mathbb{E}\Bigg(1_{\{\sup_{x\in\mathbb{Z}^d} l(n,x)\le K_n\}} \prod_{k=1}^m \varLambda_n(x_k)\Bigg)$$
$$= \sum_{\substack{i_1+\cdots+i_l=m\\ i_1,\dots,i_l\ge 1}} \frac{m!}{(i_1)!\cdots(i_l)!}\mathbb{E}\Bigg\{1_{\{\sup_{x\in\mathbb{Z}^d} l(n,x)\le K_n\}} \prod_{k=1}^l \varLambda_n(y_k)^{i_k}\Bigg\}.$$

Summarizing the above discussion,

$$\sum_{(x_1,\dots,x_m)\in A_l} \mathbb{E}\Bigg(1_{\{\sup_{x\in\mathbb{Z}^d} l(n,x)\le K_n\}} \prod_{k=1}^m \varLambda_n(x_k)\Bigg)$$
$$= \sum_{\{y_1,\dots,y_l\}\in\mathcal{C}_l} \sum_{\substack{i_1+\cdots+i_l=m\\ i_1,\dots,i_l\ge 1}} \frac{m!}{(i_1)!\cdots(i_l)!}\mathbb{E}\Bigg\{1_{\{\sup_{x\in\mathbb{Z}^d} l(n,x)\le K_n\}} \prod_{k=1}^l \varLambda_n(y_k)^{i_k}\Bigg\}.$$

Notice that the quantity

$$f(y_1,\dots,y_l) \equiv \sum_{\substack{i_1+\cdots+i_l=m\\ i_1,\dots,i_l\ge 1}} \frac{m!}{(i_1)!\cdots(i_l)!}\mathbb{E}\Bigg\{1_{\{\sup_{x\in\mathbb{Z}^d} l(n,x)\le K_n\}} \prod_{k=1}^l \varLambda_n(y_k)^{i_k}\Bigg\}$$

is invariant under the permutations over $\{y_1,\dots,y_l\}$. So we have

$$\sum_{(x_1,\dots,x_m)\in A_l} \mathbb{E}\Bigg(1_{\{\sup_{x\in\mathbb{Z}^d} l(n,x)\le K_n\}} \prod_{k=1}^m \varLambda_n(x_k)\Bigg)$$
$$= \frac{1}{l!} \sum_{(y_1,\dots,y_l)\in B_l} \sum_{\substack{i_1+\cdots+i_l=m\\ i_1,\dots,i_l\ge 1}} \frac{m!}{(i_1)!\cdots(i_l)!}\mathbb{E}\Bigg\{1_{\{\sup_{x\in\mathbb{Z}^d} l(n,x)\le K_n\}} \prod_{k=1}^l \varLambda_n(y_k)^{i_k}\Bigg\},$$

where $B_l$ is defined by (2.5). By (2.11)

$$\mathbb{E}\widetilde{H}_n^m = 2^{-m}\sum_{l=1}^m \frac{1}{l!}\sum_{\substack{i_1+\cdots+i_l=m\\ i_1,\ldots,i_l\geq 1}} \frac{m!}{(i_1)!\cdots(i_l)!} \times \sum_{(y_1,\ldots,y_l)\in B_l} \mathbb{E}\left\{1_{\{\sup_{x\in\mathbb{Z}^d} l(n,x)\leq K_n\}} \prod_{k=1}^l \Lambda_n(y_k)^{i_k}\right\}. \tag{2.12}$$

We adopt the notation "$\mathbb{E}^\omega$" for the expectation with respect to $\{\omega_k\}_{k\geq 1}$ and for each $y\in Z^d$, write $D(y)=\{1\leq k\leq n; S_k=y\}$. Then

$$\Lambda(y) = \left(\sum_{j\in D(y)}\omega_j\right)^2 - \sum_{j\in D(y)}\omega_j^2$$

and for distinct $y_1,\ldots,y_l$, the sets $D(y_1),\ldots,D(y_l)$ are disjoint. Hence, by independence,

$$\mathbb{E}^\omega \prod_{k=1}^l \Lambda_n(y_k)^{i_k} = \prod_{k=1}^l \mathbb{E}^\omega\Lambda_n(y_k)^{i_k}.$$

In particular, the above quantity is zero if any of $i_1,\ldots,i_l$ is 1. Consequently, the terms in (2.12) with $l>m/2$ are equal to zero,

$$\mathbb{E}\widetilde{H}_n = 0 \tag{2.13}$$

and for any integer $m\geq 2$,

$$\mathbb{E}\widetilde{H}_n^m = 2^{-m}\sum_{l=1}^{[2^{-1}m]} \frac{1}{l!}\sum_{\substack{i_1+\cdots+i_l=m\\ i_1,\ldots,i_l\geq 2}} \frac{m!}{(i_1)!\cdots(i_l)!} \times \sum_{(y_1,\ldots,y_l)\in B_l} \mathbb{E}\left\{1_{\{\sup_{x\in\mathbb{Z}^d} l(n,x)\leq K_n\}} \prod_{k=1}^l \mathbb{E}^\omega\Lambda_n(y_k)^{i_k}\right\}. \tag{2.14}$$

Notice that

$$\mathbb{E}^\omega\Lambda_n(y_k)^{i_k} = \mathbb{E}^\omega\left\{\left(\sum_{j=1}^{l(n,x)}\omega_j\right)^2 - \sum_{j=1}^{l(n,x)}\omega_j^2\right\}^{i_k}. \tag{2.15}$$

By Lemma 2.1 we have

$$\prod_{k=1}^l \mathbb{E}^\omega\Lambda_n(y_k)^{i_k} \leq \prod_{k=1}^l i_k! C_{i_k}^{i_k/2}\{l(n,y_k)(l(n,y_k)-1)\}^{i_k/2}$$
$$= (i_1!\cdots i_l!)(C_{i_1}^{i_1/2}\cdots C_{i_l}^{i_l/2})\prod_{k=1}^l \{l(n,y_k)(l(n,y_k)-1)\}^{i_k/2},$$

where $C_i=1$ as $i=2$ and $C_i$ is the constant $C$ given in (2.2) as $i\geq 3$. We may assume that $C\geq 1$ in the rest of the proof.

Hence,

$$\begin{aligned}\mathbb{E}\widetilde{H}_n^m &\le \frac{m!}{2^m}\sum_{l=1}^{[2^{-1}m]}\frac{1}{l!}\sum_{\substack{i_1+\cdots+i_l=m\\ i_1,\ldots,i_l\ge 2}} C_{i_1}^{i_1/2}\cdots C_{i_l}^{i_l/2}\\ &\quad\times\sum_{(y_1,\ldots,y_l)\in B_l}\mathbb{E}\biggl(1_{\{\sup_{x\in\mathbb{Z}^d} l(n,x)\le K_n\}}\prod_{k=1}^{l}\{l(n,y_k)(l(n,y_k)-1)\}^{i_k/2}\biggr)\\ &\le \frac{m!}{2^m}\sum_{l=1}^{[2^{-1}m]}\frac{1}{l!}K_n^{m-2l}\biggl\{\sum_{\substack{i_1+\cdots+i_l=m\\ i_1,\ldots,i_l\ge 2}} C_{i_1}^{i_1/2}\cdots C_{i_l}^{i_l/2}\biggr\}A_l(n)\\ &\le \frac{m!}{2^m}\sum_{l=1}^{[2^{-1}m]}\frac{1}{l!}K_n^{m-2l}2^l\biggl\{\sum_{\substack{i_1+\cdots+i_l=m\\ i_1,\ldots,i_l\ge 2}} C_{i_1}^{i_1/2}\cdots C_{i_l}^{i_l/2}\biggr\}\mathbb{E}\widetilde{Q}_n^l,\end{aligned}\tag{2.16}$$

where the last step follows from (2.6).

For each $(i_1,\ldots,i_l)$ with $i_1+\cdots+i_l=m$, write

$$k=k(i_1,\ldots,i_l)=\#\{1\le j\le l; i_j=2\}.$$

We have

$$m=i_1+\cdots+i_l\ge 2k+3(l-k)$$

which leads to $l-k\le m-2l$. Thus,

$$\begin{aligned}\sum_{\substack{i_1+\cdots+i_l=m\\ i_1,\ldots,i_l\ge 2}}(C_{i_1}^{i_1/2}\cdots C_{i_l}^{i_l/2}) &= \sum_{\substack{i_1+\cdots+i_l=m\\ i_1,\ldots,i_l\ge 2}} C^{(l-k)/2}\le C^{(m-2l)/2}\sum_{\substack{i_1+\cdots+i_l=m\\ i_1,\ldots,i_l\ge 2}} 1\\ &= C^{(m-2l)/2}\binom{m-l-1}{m-2l}.\end{aligned}$$

Hence, (2.7) follows from (2.16).

To prove (2.8), we come to (2.14) and we notice that the symmetry of $\{\omega_k\}$ implies that for any integer $l\ge 1$,

$$\mathbb{E}^{\omega}\Lambda_n(y_k)^{i_k}=2^{i_k}\mathbb{E}\biggl(\sum_{1\le j_1<j_2\le l(n,x)}\omega_{j_1}\omega_{j_2}\biggr)^{i_k}\ge 0.\tag{2.17}$$

Replacing $m$ by $2m$ in (2.4) and only keeping the term with $l=m$ on the right-hand side, we obtain

$$\begin{aligned}\mathbb{E}\widetilde{H}_n^{2m} &\ge \frac{(2m)!}{2^{2m}m!}2^{-m}\sum_{(y_1,\ldots,y_m)\in B_m}\mathbb{E}\biggl(1_{\{\sup_{x\in\mathbb{Z}^d} l(n,x)\le K_n\}}\prod_{k=1}^{m}\mathbb{E}^{\omega}\Lambda_n(y_k)^2\biggr)\\ &= \frac{(2m)!}{2^m m!}A_m(n),\end{aligned}$$

where the second step follows from (2.1) in Lemma 2.1 and (2.15).

To prove (2.9), we adopt the argument used for (2.12).

$$
\begin{aligned}
\mathbb{E}\widetilde{Q}_n^m &= 2^{-m}\sum_{x_1,\dots,x_m\in\mathbb{Z}^d}\mathbb{E}\left(1_{\{\sup_{x\in\mathbb{Z}^d}l(n,x)\le K_n\}}\prod_{k=1}^m l(n,x_k)(l(n,x_k)-1)\right)\\
&= 2^{-m}\sum_{l=1}^m\frac{1}{l!}\sum_{\substack{i_1+\cdots+i_l=m\\ i_1,\dots,i_l\ge 1}}\frac{m!}{i_1!\cdots i_l!}\\
&\quad\times\sum_{(y_1,\dots,y_l)\in B_l}\mathbb{E}\left(1_{\{\sup_{x\in\mathbb{Z}^d}l(n,x)\le K_n\}}\prod_{k=1}^l\{l(n,y_k)(l(n,y_k)-1)\}^{i_k}\right)\\
&\le\frac{m!}{2^m}\sum_{l=1}^m\frac{1}{l!}\left\{\sum_{\substack{i_1+\cdots+i_l=m\\ i_1,\dots,i_l\ge 1}}\frac{1}{i_1!\cdots i_l!}\right\}K_n^{2(m-l)}\\
&\quad\times\sum_{(y_1,\dots,y_l)\in B_l}\mathbb{E}\left(1_{\{\sup_{x\in\mathbb{Z}^d}l(n,x)\le K_n\}}\prod_{k=1}^l l(n,y_k)(l(n,y_k)-1)\right)\\
&= m!\sum_{l=1}^m\frac{1}{l!}K_n^{2(m-l)}2^{-(m-l)}A_l(n)\sum_{\substack{i_1+\cdots+i_l=m\\ i_1,\dots,i_l\ge 1}}\frac{1}{i_1!\cdots i_l!}.
\end{aligned}
$$

Finally, (2.9) follows from the following estimate:

$$
\begin{aligned}
\sum_{\substack{i_1+\cdots+i_l=m\\ i_1,\dots,i_l\ge 1}}\frac{1}{i_1!\cdots i_l!}&\le\sum_{\substack{i_1+\cdots+i_l=m\\ i_1,\dots,i_l\ge 1}}\frac{1}{(i_1-1)!\cdots(i_l-1)!}\\
&=\sum_{\substack{i_1+\cdots+i_l=m-l\\ i_1,\dots,i_l\ge 0}}\frac{1}{i_1!\cdots i_l!}=\frac{l^{m-l}}{(m-l)!}.
\end{aligned}
$$

□

**Proof of Theorem 1.1.** We start with the case $d=1$. Notice that

$$Q_n=\frac{1}{2}\left(\sum_{x\in\mathbb{Z}}l^2(n,x)-n\right).$$

By Theorem 1.2 of [6],

$$n^{-3/2}Q_n\xrightarrow{d}\frac{1}{2\sigma}\int_{-\infty}^{\infty}L^2(1,x)\,\mathrm{d}x. \tag{2.18}$$

Fix $0<\delta<1/2$ and let $K_n=n^{(1+\delta)/2}$. By the classic fact (see, for example, [19]) that

$$n^{-1/2}\sup_{x\in\mathbb{Z}}l(n,x)\xrightarrow{d}\sigma^{-1}\sup_{x\in\mathbb{R}}L(1,x) \tag{2.19}$$

we have that

$$n^{-3/2}\widetilde{Q}_n\xrightarrow{d}\frac{1}{2\sigma}\int_{-\infty}^{\infty}L^2(1,x)\,\mathrm{d}x, \tag{2.20}$$

which gives

$$\lim_{n\to\infty} n^{-3m/2}\mathbb{E}\widetilde{Q}_n^m = \frac{1}{(2\sigma)^m}\mathbb{E}\left(\int_{-\infty}^{\infty} L^2(1,x)\,\mathrm{d}x\right)^m, \quad m=1,2,\ldots. \tag{2.21}$$

Replacing $m$ by $2m+1$ in (2.7) we have

$$\begin{aligned}\mathbb{E}\widetilde{H}_n^{2m+1} &\le \frac{(2m+1)!}{2^{2m+1}}\sum_{l=1}^{m}\frac{1}{l!}n^{(1+\delta)(m-l)+(1+\delta)/2}2^l C^{(2m-2l+1)/2}\binom{2m-l}{2m-2l+1}\mathbb{E}\widetilde{Q}_n^l \\ &= \mathrm{O}\left(\sum_{l=1}^{m} n^{(1+\delta)(m-l)+(1+\delta)/2}n^{3l/2}\right) = \mathrm{o}\left(n^{3(2m+1)/4}\right), \quad n\to\infty,\end{aligned} \tag{2.22}$$

for all $m=0,1,2,\ldots$.

Replacing $m$ by $2m$ in (2.7) and by (2.21) we have

$$\begin{aligned}\mathbb{E}\widetilde{H}_n^{2m} &\le \frac{(2m)!}{2^{2m}}\sum_{l=1}^{m}\frac{1}{l!}n^{(1+\delta)(m-l)}2^l C^{m-l}\binom{2m-l-1}{2m-2l}\mathbb{E}\widetilde{Q}_n^l \\ &\sim \frac{(2m)!}{2^{2m}}\sum_{l=1}^{m}\frac{1}{l!}n^{(1+\delta)(m-l)}(2\sigma)^{-l}n^{3l/2}2^l \\ &\quad \times \mathbb{E}\left(\int_{-\infty}^{\infty} L^2(1,x)\,\mathrm{d}x\right)^l C^{m-l}\binom{2m-l-1}{2m-2l}, \quad n\to\infty.\end{aligned}$$

Clearly, the right-hand side is dominated by the term with $l=m$. Consequently,

$$\limsup_{n\to\infty} n^{-3m/2}\mathbb{E}\widetilde{H}_n^{2m} \le \frac{1}{(2\sigma)^m}\frac{(2m)!}{2^m m!}\mathbb{E}\left(\int_{-\infty}^{\infty} L^2(1,x)\,\mathrm{d}x\right)^m \tag{2.23}$$

for all $m=1,2,\ldots$. In particular, combining this with (2.8) we have

$$A_m(n) = \mathrm{O}(n^{3m/2}), \quad n\to\infty, m=1,2,\ldots. \tag{2.24}$$

On the other hand, by (2.24) and (2.21), the right-hand side of (2.9) is dominated by the term with $l=m$. Hence,

$$\liminf_{n\to\infty} n^{-3m/2}A_m(n) \ge \frac{1}{(2\sigma)^m}\mathbb{E}\left(\int_{-\infty}^{\infty} L^2(1,x)\,\mathrm{d}x\right)^m, \quad m=1,2,\ldots. \tag{2.25}$$

From (2.8),

$$\liminf_{n\to\infty} n^{-3m/2}\mathbb{E}\widetilde{H}_n^{2m} \ge \frac{1}{(2\sigma)^m}\frac{(2m)!}{2^m m!}\mathbb{E}\left(\int_{-\infty}^{\infty} L^2(1,x)\,\mathrm{d}x\right)^m. \tag{2.26}$$

In summary of (2.22), (2.23) and (2.26), and noticing that

$$\mathbb{E}U^{2m} = \frac{(2m)!}{2^m m!} \quad \text{and} \quad \mathbb{E}U^{2m+1} = 0$$

we have that for every $m=0,1,2,\ldots$,

$$\lim_{n\to\infty} n^{-3m/4}\mathbb{E}\widetilde{H}_n^m = \frac{1}{(2\sigma)^{m/2}}(\mathbb{E}U^m)\mathbb{E}\left(\int_{-\infty}^{\infty} L^2(1,x)\,\mathrm{d}x\right)^{m/2}. \tag{2.27}$$

Notice the fact that for any $\theta \in \mathbb{R}$,

$$\mathbb{E}\exp\left\{\theta\left(\int_{-\infty}^{\infty} L^2(1,x)\,\mathrm{d}x\right)^{1/2} U\right\} = \mathbb{E}\exp\left\{\frac{\theta^2}{2}\int_{-\infty}^{\infty} L^2(1,x)\,\mathrm{d}x\right\} < \infty,$$

where the last step follows from Theorem 1.1 (with $m=1$ and $p=2$) in [8]. Therefore, (2.27) implies that

$$\frac{1}{n^{3/4}}\widetilde{H}_n \xrightarrow{d} 2^{-1}\sigma^{-1/2}\left(\int_{-\infty}^{\infty} L^2(1,x)\,\mathrm{d}x\right)^{1/2} U.$$

By (2.11) and by our choice of $K_n$ we have

$$\mathbb{P}\left\{\sup_{x\in\mathbb{Z}} l(n,x) > K_n\right\} \to 0 \tag{2.28}$$

as $n\to\infty$. Thus, we have proved Theorem 1.1 in the case $d=1$.

The proof in the multi-dimensional cases is essentially the same. Instead of (2.12), we have that

$$\frac{Q_n}{n\log n} \xrightarrow{P} (2\pi)^{-1}(\det\Gamma)^{-1/2} \quad \text{as } d=2, \tag{2.29}$$

$$\frac{Q_n}{n} \xrightarrow{P} \gamma \quad \text{as } d\ge 3. \tag{2.30}$$

Indeed, (2.29) and (2.30) follow from the weak convergence of the sequences $(Q_n - \mathbb{E}Q_n)/n$ when $d=2$ (see [20]), $(Q_n-\mathbb{E}Q_n)/\sqrt{n\log n}$ when $d=3$ (see Theorem 5.1) and $(Q_n-\mathbb{E}Q_n)/\sqrt{n}$ when $d\ge 4$ (see Theorem 5.1); and from the well-known fact that

$$\mathbb{E}Q_n \sim \begin{cases} (2\pi)^{-1}(\det\Gamma)^{-1/2} n\log n, & d=2,\\ \gamma n, & d\ge 3.\end{cases} \tag{2.31}$$

In addition, it is well known [19] that

$$\sup_{x\in\mathbb{Z}^2}\frac{l(n,x)}{(\log n)^2} \quad \text{and} \quad \sup_{x\in\mathbb{Z}^d}\frac{l(n,x)}{\log n}$$

are almost surely bounded in the case $d=2$ and the case $d\ge 3$, respectively. Thus, if we define $K_n = M(\log n)^2$ as $d=2$, and $K_n = M\log n$ as $d\ge 3$. Then (2.28) holds as the constant $M>0$ is sufficiently large.

Therefore, a modification of the proof for (2.27) gives that

$$\lim_{n\to\infty}(n\log n)^{-m/2}\mathbb{E}\widetilde{H}_n^m = 2^{-m/2}(2\pi)^{-m/2}(\det\Gamma)^{-m/4}\mathbb{E}U^m \quad \text{as } d=2, \tag{2.32}$$

$$\lim_{n\to\infty} n^{-m/2}\mathbb{E}\widetilde{H}_n^m = 2^{-m/2}\gamma^{m/2}\mathbb{E}U^m \quad \text{as } d\ge 3. \tag{2.33}$$

So the multi-dimensional part of Theorem 1.1 follows from (2.32) and (2.33). □

## 3. Moderate deviations

Recall that $K_n = M(\log n)^2$ as $d=2$, where $M>0$ is a large but fixed constant. Take $K_n = (n/\log n)^{1/4}$ as $d\ge 3$. In the case $d=1$, (1.8) implies that there is a positive sequence $M_n$ such that

$$M_n\to\infty \quad \text{and} \quad M_n^2\left(\frac{b_n^7}{n}\right)^{1/4} \to 0, \quad n\to\infty. \tag{3.1}$$

So in this section we take $K_n = M_n\sqrt{nb_n}$ as $d=1$.

An important fact is that under our choice,

$$\lim_{n\to\infty}\frac{1}{b_n}\log\mathbb{P}\Big\{\sup_{x\in\mathbb{Z}^d} l(n,x)>K_n\Big\}=-\infty \tag{3.2}$$

in all dimensions. We refer to [5] for (3.2) under $d=1$, and [19] for (3.2) under $d\ge 2$.

Another important fact is that

$$\mathbb{E}\widetilde{H}_n^{2m+1}\ge 0,\quad m=0,1,\dots, \tag{3.3}$$

which follows from (2.13), (2.14) and (2.17).

We claim that Theorem 1.2 holds if we can prove that for any $\theta>0$,

$$\lim_{n\to\infty}\frac{1}{b_n}\log\mathbb{E}\exp\Big\{\pm\theta\frac{b_n^{1/4}}{n^{3/4}}\widetilde{H}_n\Big\}=\frac{\theta^4}{96\sigma^2}\qquad\text{as } d=1, \tag{3.4}$$

$$\lim_{n\to\infty}\frac{1}{b_n}\log\mathbb{E}\exp\Big\{\pm\theta\sqrt{\frac{b_n}{n\log n}}\widetilde{H}_n\Big\}=\frac{\theta^2}{4\pi\sqrt{\det\Gamma}}\qquad\text{as } d=2, \tag{3.5}$$

$$\lim_{n\to\infty}\frac{1}{b_n}\log\mathbb{E}\exp\Big\{\pm\theta\sqrt{\frac{b_n}{n}}\widetilde{H}_n\Big\}=\frac{\gamma\theta^2}{2}\qquad\text{as } d\ge 3. \tag{3.6}$$

Indeed, according to the Gärtner–Ellis theorem (Theorem 2.3.6, p. 44 in [9]), (3.4)–(3.6) imply that $\widetilde{H}_n$ satisfies the moderate deviations given in Theorem 1.2. By Theorem 4.2.13, p. 130 in [9], the moderate deviations pass from $\widetilde{H}_n$ to $H_n$ through the exponential equivalence given by

$$\limsup_{n\to\infty}\frac{1}{b_n}\log\mathbb{P}\{\widetilde{H}_n\ne H_n\}=\lim_{n\to\infty}\frac{1}{b_n}\log\mathbb{P}\Big\{\sup_{x\in\mathbb{Z}^d} l(n,x)>K_n\Big\}=-\infty,$$

where the second step follows from (3.2).

In the rest of this section, we prove (3.4), (3.5) and (3.6) in three separate parts.

*Case* $d=1$.

By (2.7) in Proposition 2.1,

$$\begin{aligned}
&\sum_{m=2}^\infty\frac{\theta^m}{m!}\Big(\frac{b_n}{n^3}\Big)^{m/4}\mathbb{E}\widetilde{H}_n^m\\
&\quad\le\sum_{m=2}^\infty\frac{\theta^m}{m!}\Big(\frac{b_n}{n^3}\Big)^{m/4}\frac{m!}{2^m}\sum_{l=1}^{[2^{-1}m]}\frac{1}{l!}K_n^{m-2l}2^lC^{(m-2l)/2}\binom{m-l-1}{m-2l}\mathbb{E}\widetilde{Q}_n^l\\
&\quad=\sum_{l=1}^\infty\frac{\theta^{2l}}{2^l l!}\Big(\frac{b_n}{n^3}\Big)^{l/2}\mathbb{E}\widetilde{Q}_n^l\sum_{m=2l}^\infty\Big(\frac{\theta}{2}\Big)^{m-2l}K_n^{m-2l}\Big(\frac{b_n}{n}\Big)^{(m-2l)/4}C^{(m-2l)/2}\binom{m-l-1}{m-2l}.
\end{aligned}$$

Notice that

$$\begin{aligned}
&\sum_{m=2l}^\infty\Big(\frac{\theta}{2}\Big)^{m-2l}K_n^{m-2l}\Big(\frac{b_n}{n^3}\Big)^{(m-2l)/4}C^{(m-2l)/2}\binom{m-l-1}{m-2l}\\
&\quad=\sum_{m=0}^\infty\Big(\frac{\theta}{2}\Big)^{m}K_n^{m}\Big(\frac{b_n}{n^3}\Big)^{m/4}C^{m/2}\binom{m+l-1}{m}\\
&\quad=\Big(1-\frac{\sqrt{C}\theta K_n b_n^{1/4}}{2n^{3/4}}\Big)^{-(l-1)},
\end{aligned}$$

where the last step follows from the Taylor expansion:

$$(1-x)^{-(l-1)} = \sum_{m=0}^{\infty} \binom{m+l-1}{m} x^m, \quad |x|<1. \tag{3.7}$$

Combining the above estimate with (2.13) gives

$$\mathbb{E}\exp\left\{\theta \frac{b_n^{1/4}}{n^{3/4}} \widetilde{H}_n\right\} \le \mathbb{E}\exp\left\{\frac{\theta^2}{2}\frac{b_n^{1/2}}{n^{3/2}}\left(1-\frac{\sqrt{C}\theta K_n b_n^{1/4}}{2n^{3/4}}\right)^{-1}\widetilde{Q}_n\right\}.$$

In view of (2.13), by the Taylor expansion one can easily see that

$$\mathbb{E}\exp\left\{-\theta \frac{b_n^{1/4}}{n^{3/4}} \widetilde{H}_n\right\} \le \mathbb{E}\exp\left\{\theta \frac{b_n^{1/4}}{n^{3/4}} \widetilde{H}_n\right\}.$$

So we have

$$\mathbb{E}\exp\left\{\pm\theta \frac{b_n^{1/4}}{n^{3/4}} \widetilde{H}_n\right\} \le \mathbb{E}\exp\left\{\frac{\theta^2}{2}\frac{b_n^{1/2}}{n^{3/2}}\left(1-\frac{\sqrt{C}\theta K_n b_n^{1/4}}{2n^{3/4}}\right)^{-1}\widetilde{Q}_n\right\}. \tag{3.8}$$

Notice that

$$Q_n \le \frac{1}{2}\sum_{x\in\mathbb{Z}} l^2(n,x) \le \frac{1}{2} n \sup_{x\in\mathbb{Z}} l(n,x). \tag{3.9}$$

For any $\lambda > 0$,

$$\mathbb{E}\exp\left\{\lambda \frac{b_n^{1/2}}{n^{3/2}} Q_n\right\} \le \mathbb{E}\exp\left\{\frac{\lambda}{2}\sqrt{\frac{b_n}{n}} \sup_{x\in\mathbb{Z}} l(n,x)\right\}.$$

By the fact (see Lemmas 11 and 12 in [15]) that

$$\limsup_{n\to\infty} \frac{1}{b_n} \log \mathbb{E}\exp\left\{\frac{\lambda}{2}\sqrt{\frac{b_n}{n}} \sup_{x\in\mathbb{Z}} l(n,x)\right\} < \infty,$$

we have that for any $\lambda > 0$

$$\limsup_{n\to\infty} \frac{1}{b_n} \log \mathbb{E}\exp\left\{\lambda \frac{b_n^{1/2}}{n^{3/2}} Q_n\right\} < \infty. \tag{3.10}$$

Recall (Theorem 1.3 in [8], with $m=1$ and $p=2$) that for any $\lambda > 0$

$$\lim_{n\to\infty} \frac{1}{b_n} \log \mathbb{P}\{Q_n \ge \lambda n^{3/2} b_n^{1/2}\} = -6\sigma^2\lambda^2. \tag{3.11}$$

According to Varadhan's integral lemma (Theorem 4.3.1, p. 137 in [9]), (3.10) and (3.11) imply that for any $\lambda > 0$,

$$\lim_{n\to\infty} \frac{1}{b_n} \log \mathbb{E}\exp\left\{\lambda \frac{b_n^{1/2}}{n^{3/2}} Q_n\right\} = \sup_{y>0}\{y\lambda - 6\sigma^2\lambda^2\} = \frac{\lambda^2}{24\sigma^2}. \tag{3.12}$$

This, together with (3.8), gives the desired upper bound for (3.4):

$$\limsup_{n\to\infty} \frac{1}{b_n} \log \mathbb{E}\exp\left\{\pm\theta \frac{b_n^{1/4}}{n^{3/4}} \widetilde{H}_n\right\} \le \frac{\theta^4}{96\sigma^2}. \tag{3.13}$$

On the other hand, by (2.8)

$$\sum_{m=1}^{\infty}\frac{\theta^{2m}}{(2m)!}\left(\frac{b_n}{n^3}\right)^{m/2}\mathbb{E}\widetilde{H}_n^{2m}\geq\sum_{m=1}^{\infty}\frac{1}{m!}\frac{\theta^{2m}}{2^m}\left(\frac{b_n}{n^3}\right)^{m/2}A_m(n). \tag{3.14}$$

Write

$$\bar{\theta}=\theta\exp\left\{-\frac{\theta^2K_n^2}{8}\sqrt{\frac{b_n}{n^3}}\right\}.$$

By (2.9),

$$\begin{aligned}
&\sum_{m=1}^{\infty}\frac{1}{m!}\frac{\bar{\theta}^{2m}}{2^m}\left(\frac{b_n}{n^3}\right)^{m/2}\mathbb{E}\widetilde{Q}_n^m\\
&\quad\leq\sum_{m=1}^{\infty}\frac{\bar{\theta}^{2m}}{2^m}\left(\frac{b_n}{n^3}\right)^{m/2}\sum_{l=1}^{m}\binom{m}{l}\left(\frac{lK_n^2}{2}\right)^{m-l}A_l(n)\\
&\quad=\sum_{l=1}^{\infty}\frac{1}{l!}\left(\frac{\bar{\theta}^2}{2}\right)^{l}\left(\frac{b_n}{n^3}\right)^{l/2}A_l(n)\sum_{m=l}^{\infty}\frac{1}{(m-l)!}\left(\frac{\bar{\theta}^2}{2}\right)^{m-l}\left(\frac{b_n}{n^3}\right)^{(m-l)/2}\left(\frac{lK_n^2}{2}\right)^{m-l}\\
&\quad=\sum_{l=1}^{\infty}\frac{1}{l!}\left(\frac{\bar{\theta}^2}{2}\right)^{l}\left(\frac{b_n}{n^3}\right)^{l/2}A_l(n)\exp\left\{l\frac{\bar{\theta}^2K_n^2}{4}\sqrt{\frac{b_n}{n^3}}\right\}\\
&\quad\leq\sum_{l=1}^{\infty}\frac{1}{l!}\left(\frac{1}{2}\right)^{l}\left(\bar{\theta}\exp\left\{\frac{\theta^2K_n^2}{8}\sqrt{\frac{b_n}{n^3}}\right\}\right)^{2l}\left(\frac{b_n}{n^3}\right)^{l/2}A_l(n)\\
&\quad=\sum_{l=1}^{\infty}\frac{1}{l!}\frac{\theta^{2l}}{2^l}\left(\frac{b_n}{n^3}\right)^{l/2}A_l(n),
\end{aligned} \tag{3.15}$$

where the second inequality follows from the fact that $\bar{\theta}\leq\theta$.

Combining (3.14) and (3.15),

$$\begin{aligned}
\sum_{m=0}^{\infty}\frac{\theta^{2m}}{(2m)!}\left(\frac{b_n}{n^3}\right)^{m/2}\mathbb{E}\widetilde{H}_n^{2m}&\geq\mathbb{E}\exp\left\{\frac{\bar{\theta}^2}{2}\frac{b_n^{1/2}}{n^{3/2}}\widetilde{Q}_n\right\}\\
&=(1+\mathrm{o}(1))\mathbb{E}\exp\left\{\frac{\theta^2}{2}\frac{b_n^{1/2}}{n^{3/2}}\widetilde{Q}_n\right\},
\end{aligned} \tag{3.16}$$

where the second step follows from the estimate (notice that (3.9) implies that $\widetilde{Q}_n\leq nK_n/2$)

$$\begin{aligned}
\frac{\bar{\theta}^2}{2}\frac{b_n^{1/2}}{n^{3/2}}\widetilde{Q}_n&=\frac{\theta^2}{2}\frac{b_n^{1/2}}{n^{3/2}}\widetilde{Q}_n-\frac{\theta^2}{2}\left[1-\exp\left\{\frac{\theta^2K_n^2}{4}\sqrt{\frac{b_n}{n^3}}\right\}\right]\frac{b_n^{1/2}}{n^{3/2}}\widetilde{Q}_n\\
&\geq\frac{\theta^2}{2}\frac{b_n^{1/2}}{n^{3/2}}\widetilde{Q}_n-\frac{\theta^2}{2}\frac{\theta^2K_n^2}{4}\sqrt{\frac{b_n}{n^3}}\frac{b_n^{1/2}}{n^{3/2}}\frac{nK_n}{2}=\frac{\theta^2}{2}\frac{b_n^{1/2}}{n^{3/2}}\widetilde{Q}_n-\mathrm{o}(1),\quad n\to\infty.
\end{aligned} \tag{3.17}$$

In view of (3.3), by (3.16) we conclude that

$$\mathbb{E}\exp\left\{\theta\frac{b_n^{1/4}}{n^{3/4}}\widetilde{H}_n\right\}\geq(1+\mathrm{o}(1))\mathbb{E}\exp\left\{\frac{\theta^2}{2}\frac{b_n^{1/2}}{n^{3/2}}\widetilde{Q}_n\right\}. \tag{3.18}$$

Getting the lower bound for the negative coefficients $-\theta$ is harder. To do this we need to control the terms with odd powers. Replacing $m$ by $2m+1$ in (2.7) we obtain

$$\begin{aligned}
&\sum_{m=1}^{\infty}\frac{\theta^{2m+1}}{(2m+1)!}\left(\frac{b_n}{n^3}\right)^{(2m+1)/4}\mathbb{E}\widetilde{H}_n^{2m+1}\\
&\quad\le\sum_{m=1}^{\infty}\frac{\theta^{2m+1}}{2^{2m+1}}\left(\frac{b_n}{n^3}\right)^{(2m+1)/4}\sum_{l=1}^{m}\frac{1}{l!}K_n^{2(m-l)+1}2^lC^{(2m-2l+1)/2}\binom{2m-l}{2m-2l+1}\mathbb{E}\widetilde{Q}_n^l\\
&\quad=\frac{\theta\sqrt{C}}{2}K_n\left(\frac{b_n}{n^3}\right)^{1/4}\sum_{l=1}^{\infty}\frac{1}{l!}\frac{\theta^{2l}}{2^l}\left(\frac{b_n}{n^3}\right)^{l/2}\mathbb{E}\widetilde{Q}_n^l\\
&\qquad\times\sum_{m=l}^{\infty}\left(\frac{\theta}{2}\right)^{2(m-l)}C^{m-l}K_n^{2(m-l)}\left(\frac{b_n}{n^3}\right)^{(m-l)/2}\binom{2m-l}{2m-2l+1}.
\end{aligned}\tag{3.19}$$

Noticing that

$$\binom{2m-l}{2m-2l+1}=\frac{2m-l}{2m-2l+1}\binom{2m-l-1}{2m-2l}\le l\binom{2m-l-1}{2m-2l}$$

we have

$$\begin{aligned}
&\sum_{m=l}^{\infty}\left(\frac{\theta}{2}\right)^{2(m-l)}C^{m-l}K_n^{2(m-l)}\left(\frac{b_n}{n^3}\right)^{(m-l)/2}\binom{2m-l-1}{2m-2l}\\
&\quad=l\sum_{m=0}^{\infty}\left(\frac{\theta}{2}\right)^{2m}K_n^{2m}\left(\frac{b_n}{n^3}\right)^{m/2}C^m\binom{2m+l-1}{2m}\\
&\quad\le l\sum_{m=0}^{\infty}\left(\frac{\theta}{2}\right)^{m}K_n^{m}\left(\frac{b_n}{n^3}\right)^{m/4}C^m\binom{m+l-1}{m}\\
&\quad=l\left(1-\frac{\sqrt{C}\theta K_nb_n^{1/4}}{2n^{3/4}}\right)^{-(l-1)},
\end{aligned}$$

where the last step follows from (3.7).

By (3.17), therefore,

$$\begin{aligned}
&\sum_{m=1}^{\infty}\frac{\theta^{2m+1}}{(2m+1)!}\left(\frac{b_n}{n^3}\right)^{(2m+1)/4}\mathbb{E}\widetilde{H}_n^{2m+1}\\
&\quad\le\sqrt{C}\frac{\theta^3}{4}K_n\left(\frac{b_n}{n^3}\right)^{3/4}\sum_{l=1}^{\infty}\frac{1}{(l-1)!}\left(\frac{\theta^2}{2}\right)^{l-1}\left(\frac{b_n}{n^3}\right)^{(l-1)/2}\\
&\qquad\times\left(1-\frac{\sqrt{C}\theta K_nb_n^{1/4}}{2n^{3/4}}\right)^{-(l-1)}\mathbb{E}\widetilde{Q}_n^l\\
&\quad\le\sqrt{C}\frac{\theta^3}{4}K_n\left(\frac{b_n}{n^3}\right)^{3/4}\mathbb{E}\left(\widetilde{Q}_n\exp\left\{\frac{\theta^2}{2}\frac{b_n^{1/2}}{n^{3/2}}\left(1-\frac{\sqrt{C}\theta K_nb_n^{1/4}}{2n^{3/4}}\right)^{-1}\widetilde{Q}_n\right\}\right).
\end{aligned}$$

By the estimate $\widetilde{Q}_n\le\frac{1}{2}nK_n$ and by the assumption (1.8) we have that

$$K_n\left(\frac{b_n}{n^3}\right)^{3/4}Q_n\le\frac{M^2}{2}\frac{b_n^{7/4}}{n^{1/4}}\to 0,\quad n\to\infty.$$

In view of (2.13), we have proved that

$$\sum_{m=0}^{\infty}\frac{\theta^{2m+1}}{(2m+1)!}\left(\frac{b_n}{n^3}\right)^{(2m+1)/4}\mathbb{E}\widetilde{H}_n^{2m+1}=\mathrm{o}\left(\mathbb{E}\exp\left\{\frac{\theta^2}{2}\frac{b_n^{1/2}}{n^{3/2}}\widetilde{Q}_n\right\}\right),\quad n\to\infty. \tag{3.20}$$

This, together with (3.16), yields

$$\mathbb{E}\exp\left\{-\theta\frac{b_n^{1/4}}{n^{3/4}}\widetilde{H}_n\right\}\geq(1+\mathrm{o}(1))\mathbb{E}\exp\left\{\frac{\theta^2}{2}\frac{b_n^{1/2}}{n^{3/2}}\widetilde{Q}_n\right\},\quad n\to\infty.$$

Combining this with (3.18),

$$\mathbb{E}\exp\left\{\pm\theta\frac{b_n^{1/4}}{n^{3/4}}\widetilde{H}_n\right\}\geq(1+\mathrm{o}(1))\mathbb{E}\exp\left\{\frac{\theta^2}{2}\frac{b_n^{1/2}}{n^{3/2}}\widetilde{Q}_n\right\},\quad n\to\infty. \tag{3.21}$$

To estimate the right-hand side of (3.21), notice that

$$\begin{aligned}\mathbb{E}\exp\left\{\frac{\theta^2}{2}\frac{b_n^{1/2}}{n^{3/2}}\widetilde{Q}_n\right\}&\geq\mathbb{E}\left(\exp\left\{\frac{\theta^2}{2}\frac{b_n^{1/2}}{n^{3/2}}Q_n\right\}1_{\{\sup_{x\in\mathbb{Z}}l(n,x)\leq K_n\}}\right)\\&=\mathbb{E}\exp\left\{\frac{\theta^2}{2}\frac{b_n^{1/2}}{n^{3/2}}Q_n\right\}-\mathbb{E}\left(\exp\left\{\frac{\theta^2}{2}\frac{b_n^{1/2}}{n^{3/2}}Q_n\right\}1_{\{\sup_{x\in\mathbb{Z}}l(n,x)>K_n\}}\right).\end{aligned}$$

Consequently, by (3.12)

$$\begin{aligned}&\max\left\{\liminf_{n\to\infty}\frac{1}{b_n}\log\mathbb{E}\exp\left\{\frac{\theta^2}{4}\frac{b_n^{1/2}}{n^{3/2}}\widetilde{Q}_n\right\},\limsup_{n\to\infty}\frac{1}{b_n}\log\mathbb{E}\left(\exp\left\{\frac{\theta^2}{4}\frac{b_n^{1/2}}{n^{3/2}}Q_n\right\}1_{\{\sup_{x\in\mathbb{Z}}l(n,x)>K_n\}}\right)\right\}\\&\quad\geq\frac{\theta^4}{96\sigma^2}.\end{aligned} \tag{3.22}$$

By the Cauchy–Schwarz inequality,

$$\begin{aligned}&\mathbb{E}\left(\exp\left\{\frac{\theta^2}{2}\frac{b_n^{1/2}}{n^{3/2}}Q_n\right\}1_{\{\sup_{x\in\mathbb{Z}}l(n,x)>K_n\}}\right)\\&\quad\leq\left(\mathbb{E}\exp\left\{\theta^2\frac{b_n^{1/2}}{n^{3/2}}Q_n\right\}\right)^{1/2}\left(\mathbb{P}\left\{\sup_{x\in\mathbb{Z}}l(n,x)>K_n\right\}\right)^{1/2}.\end{aligned}$$

Hence, (3.2) and (3.12) imply that

$$\limsup_{n\to\infty}\frac{1}{b_n}\log\mathbb{E}\left(\exp\left\{\frac{\theta^2}{2}\frac{b_n^{1/2}}{n^{3/2}}Q_n\right\}1_{\{\sup_{x\in\mathbb{Z}}l(n,x)>K_n\}}\right)=-\infty.$$

In view of (3.22),

$$\liminf_{n\to\infty}\frac{1}{b_n}\log\mathbb{E}\exp\left\{\frac{\theta^2}{4}\frac{b_n^{1/2}}{n^{3/2}}\widetilde{Q}_n\right\}\geq\frac{\theta^4}{96\sigma^2}. \tag{3.23}$$

Combining (3.21) and (3.23) gives the desired lower bound

$$\liminf_{n\to\infty}\frac{1}{b_n}\log\mathbb{E}\exp\left\{\pm\theta\frac{b_n^{1/4}}{n^{3/4}}\widetilde{H}_n\right\}\geq\frac{\theta^4}{96\sigma^2}. \tag{3.24}$$

Therefore, (3.4) follows from (3.13) and (3.24).

*Case* $d = 2$.

Similar to (3.8) and (3.21), respectively,

$$\mathbb{E}\exp\left\{\pm\theta\sqrt{\frac{b_n}{n\log n}}\widetilde{H}_n\right\} \le \mathbb{E}\exp\left\{\frac{\theta^2}{2}\left(1-\frac{\sqrt{C}\theta K_n}{2}\sqrt{\frac{b_n}{n\log n}}\right)^{-1}\frac{b_n}{n\log n}\widetilde{Q}_n\right\} \tag{3.25}$$

and

$$\mathbb{E}\exp\left\{\pm\theta\sqrt{\frac{b_n}{n\log n}}\widetilde{H}_n\right\} \ge (1+\mathrm{o}(1))\mathbb{E}\exp\left\{\frac{\theta^2}{2}\frac{b_n}{n\log n}\widetilde{Q}_n\right\}, \quad n\to\infty. \tag{3.26}$$

Applying Jensen's inequality on the right-hand side of (3.26),

$$\mathbb{E}\exp\left\{\theta\sqrt{\frac{b_n}{n\log n}}\widetilde{H}_n\right\} \ge (1+\mathrm{o}(1))\exp\left\{\frac{\theta^2}{2}\frac{b_n}{n\log n}\mathbb{E}\widetilde{Q}_n\right\}, \quad \theta\in\mathbb{R}. \tag{3.27}$$

By the fact (implied by (3.2)) that $\mathbb{P}\{\sup_{x\in\mathbb{Z}^2} l(n,x) > K_n\} \longrightarrow 0$, $n\to\infty$, we have that

$$\mathbb{E}\widetilde{Q}_n \sim \mathbb{E}Q_n \sim (2\pi)^{-1}\det\Gamma^{-1/2} n\log n, \quad n\to\infty,$$

where the second step follows from (2.31).

Consequently,

$$\liminf_{n\to\infty}\frac{1}{b_n}\log\mathbb{E}\exp\left\{\pm\theta\sqrt{\frac{b_n}{n\log n}}\widetilde{H}_n\right\} \ge \frac{\theta^2}{4\pi\sqrt{\det\Gamma}}. \tag{3.28}$$

On the other hand, recall the fact (Lemma 2.3 in [3]) that

$$\mathbb{E}\exp\left\{\frac{\lambda}{n}|Q_n - \mathbb{E}Q_n|\right\} < \infty$$

for some $\lambda > 0$. By the assumption (1.10) we have

$$\begin{aligned}\mathbb{E}\exp\left\{\frac{\theta^2}{2}\frac{b_n}{n\log n}Q_n\right\} &\le \exp\left\{\frac{\theta^2}{2}\frac{b_n}{n\log n}\mathbb{E}Q_n\right\}\mathbb{E}\exp\left\{\frac{\theta^2}{2}\frac{b_n}{n\log n}|Q_n - \mathbb{E}Q_n|\right\}\\ &= \mathrm{O}\left(\exp\left\{\frac{\theta^2}{2}\frac{b_n}{n\log n}\mathbb{E}Q_n\right\}\right), \quad n\to\infty.\end{aligned} \tag{3.29}$$

Combining this with (3.25) gives

$$\limsup_{n\to\infty}\frac{1}{b_n}\log\mathbb{E}\exp\left\{\pm\theta\sqrt{\frac{b_n}{n\log n}}\widetilde{H}_n\right\} \le \frac{\theta^2}{4\pi\sqrt{\det\Gamma}}. \tag{3.30}$$

Thus, (3.5) follows from (3.28) and (3.30).

*Case* $d \ge 3$.

The treatment in the case $d \ge 3$ is almost same as the one given in the case $d = 2$, except that here we use

$$\mathbb{E}\exp\left\{\frac{\theta^2}{2}\frac{b_n}{n}\widetilde{Q}_n\right\} = \mathrm{O}\left(\exp\left\{\frac{\theta^2}{2}\frac{b_n}{n}\mathbb{E}Q_n\right\}\right), \quad n\to\infty, \tag{3.31}$$

instead of (3.29).

We end this section with the proof of (3.31). Notice that

$$\mathbb{E}\exp\left\{\frac{\theta^2}{2}\frac{b_n}{n}\widetilde{Q}_n\right\} \le 1 + \mathbb{E}\left(\exp\left\{\frac{\theta^2}{2}\frac{b_n}{n}Q_n\right\}1_{\{\sup_{x\in\mathbb{Z}^d} l(n,x)\le K_n\}}\right)$$
$$\le 1 + \mathbb{E}\exp\left\{\frac{\theta^2}{2}\frac{b_n}{n}\mathbb{E}\widetilde{Q}_n\right\}\mathbb{E}\left(\exp\left\{\frac{\theta^2}{2}\frac{b_n}{n}|Q_n - \mathbb{E}Q_n|\right\}1_{\{\sup_{x\in\mathbb{Z}^d} l(n,x)\le K_n\}}\right).$$

Therefore, we need only to prove that

$$\sup_n \mathbb{E}\left(\exp\left\{\frac{\theta^2}{2}\frac{b_n}{n}|Q_n - \mathbb{E}Q_n|\right\}1_{\{\sup_{x\in\mathbb{Z}^d} l(n,x)\le K_n\}}\right) < \infty. \tag{3.32}$$

By the fact that $Q_n \le 2^{-1}nK_n$ on the event $\{\sup_{x\in\mathbb{Z}^d} l(n,x) \le K_n\}$, we have $|Q_n - \mathbb{E}Q_n| \le 2^{-1}nK_n$. Consequently,

$$\mathbb{E}\left(\exp\left\{\frac{\theta^2}{2}\frac{b_n}{n}|Q_n - \mathbb{E}Q_n|\right\}1_{\{\sup_{x\in\mathbb{Z}^d} l(n,x)\le K_n\}}\right)$$
$$\le \mathbb{E}\exp\left\{\frac{\theta^2}{2}\frac{b_n}{n}(2^{-1}nK_n)^{1/3}|Q_n - \mathbb{E}Q_n|^{2/3}\right\}.$$

Finally, (3.32) follows from Theorem 5.2 and our assumptions on $\{b_n\}$ given in Theorem 1.2.

## 4. Laws of the iterated logarithm

The following Lévy type inequality is needed in our proof of the upper bounds in Theorem 1.3.

**Lemma 4.1.** *For any $s, t > 0$ and integer $n \ge 2$,*

$$\min_{1\le k\le n} \mathbb{P}\{|H_k| \le s\}\mathbb{P}\left\{\max_{1\le l\le n}|H_l| \ge s + t\right\} \le 2\mathbb{P}\{|H_n| \ge t\}. \tag{4.1}$$

**Proof.** Write

$$\tau = \inf\{l \ge 1; |H_l| \ge s + t\}$$

and notice that for each $1 \le l \le n$,

$$\sum_{l+1\le j<k\le n} \omega_j\omega_k 1_{\{S_j=S_k\}} \overset{d}{=} H_{n-l},$$

where we follow the convention that both sides are zero if $l = n-1$ or $n = l$. Thus

$$\min_{1\le k\le n} \mathbb{P}\{|H_k| \le s\}\mathbb{P}\left\{\max_{1\le l\le n}|H_l| \ge s+t\right\}$$
$$= \sum_{l=1}^n \min_{1\le k\le n} \mathbb{P}\{|H_k| \le s\}\mathbb{P}\{\tau = l\}$$
$$\le \sum_{l=1}^n \mathbb{P}\left\{\left|\sum_{l+1\le j<k\le n}\omega_j\omega_k 1_{\{S_j=S_k\}}\right| \le s\right\}\mathbb{P}\{\tau = l\}$$
$$= \sum_{l=1}^n \mathbb{P}\left\{\tau = l, \left|\sum_{l+1\le j<k\le n}\omega_j\omega_k 1_{\{S_j=S_k\}}\right| \le s\right\}, \tag{4.2}$$

where the last step follows from independence between $\{\tau = l\}$ and

$$\sum_{l+1\le j<k\le n} \omega_j\omega_k 1_{\{S_j=S_k\}}.$$

For each $1 \le l \le n$, write

$$H_n^{(l)} = \sum_{1\le j<k\le l} \omega_j\omega_k 1_{\{S_j=S_k\}} + \sum_{l+1\le j<k\le n} \omega_j\omega_k 1_{\{S_j=S_k\}} - \sum_{j=1}^{l}\sum_{k=l+1}^{n} \omega_j\omega_k 1_{\{S_j=S_k\}}.$$

Notice that

$$\begin{aligned}
&\left\{\tau = l, \left|\sum_{l+1\le j<k\le n} \omega_j\omega_k 1_{\{S_j=S_k\}}\right| \le s\right\} \\
&\quad \subset \left\{\tau = l, \left|\sum_{1\le j<k\le l} \omega_j\omega_k 1_{\{S_j=S_k\}} + \sum_{l+1\le j<k\le n} \omega_j\omega_k 1_{\{S_j=S_k\}}\right| \ge t\right\} \\
&\quad \subset \mathbb{P}\{\tau = l, |H_n| + |H_n^{(l)}| \ge 2t\}.
\end{aligned}$$

Hence,

$$\begin{aligned}
&\mathbb{P}\left\{\tau = l, \left|\sum_{l+1\le j<k\le n} \omega_j\omega_k 1_{\{S_j=S_k\}}\right| \le s\right\} \\
&\quad \le \mathbb{P}\{\tau = l, |H_n| \ge t\} + \mathbb{P}\{\tau = l, |H_n^{(l)}| \ge t\}.
\end{aligned} \tag{4.3}$$

We now claim that

$$\mathbb{P}\{\tau = l, |H_n^{(l)}| \ge t\} = \mathbb{P}\{\tau = l, |H_n| \ge t\}. \tag{4.4}$$

Indeed, (4.4) follows from the fact that the random vectors

$$(\omega_1, \dots, \omega_n) \quad \text{and} \quad (\omega_1, \dots, \omega_l, -\omega_{l+1}, \dots, -\omega_n)$$

have the same distribution, and that replacing the first vector by the second does not change the event $\{\tau = l\}$ but changes $H_n$ into $H_n^{(l)}$.

Finally, by (4.2), (4.3) and (4.4),

$$\begin{aligned}
&\min_{1\le k\le n} \mathbb{P}\{|H_k| \le s\}\mathbb{P}\left\{\max_{1\le l\le n} |H_l| \ge s + t\right\} \\
&\quad \le 2\sum_{l=1}^{n} \mathbb{P}\{\tau = l, |H_n| \ge t\} \le 2\mathbb{P}\{|H_n| \ge t\}.
\end{aligned}$$

□

**Proof of Theorem 1.3.** Due to similarity we only consider the case $d = 1$. To prove the upper bound in (1.13), it suffices to show

$$\limsup_{n\to\infty} \frac{|H_n|}{(n \log\log n)^{3/4}} \le \frac{2^{3/4}}{3}\sigma^{-1/2} \quad \text{a.s.} \tag{4.5}$$

Let $\theta > 0$ and

$$\lambda_1 > \frac{2^{3/4}}{3}\sigma^{-1/2}$$

be fixed but arbitrary. Write $n_k = [\theta^k]$ for $k = 1, 2, \ldots$. Take $\varepsilon > 0$ small enough so

$$\lambda_1 - \varepsilon > \frac{2^{3/4}}{3}\sigma^{-1/2}.$$

By Theorem 1.1,

$$\min_{1\le m\le n} \mathbb{P}\{|H_m| \le \varepsilon(n_k \log\log n_k)^{3/4}\} \ge \frac{1}{2}$$

as $k$ is sufficiently large. By Lemma 4.1, therefore,

$$\mathbb{P}\Big\{\max_{1\le l\le n_k} |H_l| \ge \lambda_1(n_k \log\log n_k)^{3/4}\Big\} \le 4\mathbb{P}\{|H_{n_k}| \ge (\lambda_1 - \varepsilon)(n_k \log\log n_k)^{3/4}\}.$$

By (1.7) in Theorem 1.2 (with $b_n = \log\log n$),

$$\sum_k \mathbb{P}\Big\{\max_{1\le l\le n_k} |H_l| \ge \lambda_1(n_k \log\log n_k)^{3/4}\Big\} < \infty.$$

By the Borel–Cantelli lemma,

$$\limsup_{k\to\infty} \frac{1}{(n_k \log\log n_k)^{3/4}} \max_{1\le l\le n_k} |H_l| \le \lambda_1 \quad \text{a.s.}$$

For any large integer $n$, if $n_k \le n \le n_{k+1}$, then

$$\frac{|H_n|}{(n\log\log n)^{3/4}} \le (\theta^{3/4} + \mathrm{o}(1))\frac{1}{(n_{k+1}\log\log n_{k+1})^{3/4}} \max_{1\le l\le n_{k+1}} |H_l|.$$

So we have

$$\limsup_{n\to\infty} \frac{|H_n|}{(n\log\log n)^{3/4}} \le \theta^{3/4}\lambda_1 \quad \text{a.s.}$$

Letting $\theta \to 1^+$ and $\lambda_1 \to 2^{3/4}3^{-1}\sigma^{-1/2}$ gives (4.5).

We only prove the lower bound for $H_n$:

$$\limsup_{n\to\infty} \frac{H_n}{(n\log\log n)^{3/4}} \ge \frac{2^{3/4}}{3}\sigma^{-1/2} \quad \text{a.s.} \tag{4.6}$$

as the proof of the lower bound for $-H_n$ is analogous.

Let $n_k$ be defined as above (but with large constant $\theta > 0$) and let the constant $\lambda_2$ satisfying

$$\lambda_2 < \frac{2^{3/4}}{3}\sigma^{-1/2}.$$

Let $\varepsilon > 0$ be small enough so

$$\lambda_2 + \varepsilon < \frac{2^{3/4}}{3}\sigma^{-1/2}.$$

Notice that

$$\sum_{n_k+1\le i<j\le n_{k+1}} \omega_i\omega_j 1_{\{S_i=S_j\}} \overset{d}{=} H_{n_{k+1}-n_k}.$$

As $\theta > 0$ and $k$ are sufficiently large,

$$\mathbb{P}\Big\{\sum_{n_k+1\le i<j\le n_{k+1}} \omega_i\omega_j 1_{\{S_i=S_j\}} \ge \lambda_2 (n_{k+1}\log\log n_{k+1})^{3/4}\Big\}$$
$$\ge \mathbb{P}\{H_{n_{k+1}-n_k} \ge (\lambda_2+\varepsilon)((n_{k+1}-n_k)\log\log(n_{k+1}-n_k))^{3/4}\}.$$

By (1.7) in Theorem 1.2 again,

$$\sum_k \mathbb{P}\Big\{\sum_{n_k+1\le i<j\le n_{k+1}} \omega_i\omega_j 1_{\{S_i=S_j\}} \ge \lambda_2 (n_{k+1}\log\log n_{k+1})^{3/4}\Big\} = \infty.$$

Notice that

$$\sum_{n_k+1\le i<j\le n_{k+1}} \omega_i\omega_j 1_{\{S_i=S_j\}}, \quad k=1,2,\dots,$$

is an independent sequence. By the Borel–Cantelli lemma,

$$\limsup_{k\to\infty} \frac{1}{(n_{k+1}\log\log n_{k+1})^{3/4}} \sum_{n_k+1\le i<j\le n_{k+1}} \omega_i\omega_j 1_{\{S_i=S_j\}} \ge \lambda_2 \quad \text{a.s.}$$

In addition, (4.5) implies that

$$\limsup_{k\to\infty} \frac{1}{(n_{k+1}\log\log n_{k+1})^{3/4}} |H_{n_k}| \le \theta^{-3/4}\frac{2^{3/4}}{3}\sigma^{-1/2} \quad \text{a.s.}$$

Consequently,

$$\limsup_{k\to\infty} \frac{1}{(n_{k+1}\log\log n_{k+1})^{3/4}} \Big\{H_{n_k} + \sum_{n_k+1\le i<j\le n_{k+1}} \omega_i\omega_j 1_{\{S_i=S_j\}}\Big\}$$
$$\ge \lambda_2 - \theta^{-3/4}\frac{2^{3/4}}{3} \quad \text{a.s.} \tag{4.7}$$

Recall the notation

$$H_{n_{k+1}}^{(n_k)} = H_{n_k} + \sum_{n_k+1\le i<j\le n_{k+1}} \omega_i\omega_j 1_{\{S_i=S_j\}} - \sum_{i=1}^{n_k}\sum_{j=n_k+1}^{n_{k+1}} 1_{\{S_i=S_j\}}.$$

We have

$$H_{n_k} + \sum_{n_k+1\le i<j\le n_{k+1}} \omega_i\omega_j 1_{\{S_i=S_j\}} = \frac{H_{n_{k+1}} + H_{n_{k+1}}^{(n_k)}}{2}.$$

Therefore, by (4.7)

$$\limsup_{k\to\infty} \frac{H_{n_{k+1}}}{(n_{k+1}\log\log n_{k+1})^{3/4}} + \limsup_{k\to\infty} \frac{H_{n_{k+1}}^{(n_k)}}{(n_{k+1}\log\log n_{k+1})^{3/4}}$$
$$\ge 2\Big(\lambda_2 - \theta^{-3/4}\frac{2^{3/4}}{3}\sigma^{-1/2}\Big) \quad \text{a.s.} \tag{4.8}$$

On the other hand, notice that for each $k$,

$$H_{n_{k+1}}^{(n_k)} \stackrel{d}{=} H_{n_{k+1}}.$$

By Theorem 1.2 (with $b_n = \log\log n$)

$$\sum_k \mathbb{P}\{H_{n_{k+1}}^{(n_k)} \ge \lambda(n_{k+1}\log\log n_{k+1})^{3/4}\}$$
$$= \sum_k \mathbb{P}\{H_{n_{k+1}} \ge \lambda(n_{k+1}\log\log n_{k+1})^{3/4}\} < \infty \quad \forall\lambda > \frac{2^{3/4}}{3}\sigma^{-1/2}.$$

By the Borel–Cantelli lemma,

$$\limsup_{k\to\infty} \frac{H_{n_{k+1}}^{(n_k)}}{(n_{k+1}\log\log n_{k+1})^{3/4}} \le \frac{2^{3/4}}{3}\sigma^{-1/2} \quad \text{a.s.}$$

Combining this with (4.8) yields

$$\limsup_{k\to\infty} \frac{H_{n_{k+1}}}{(n_{k+1}\log\log n_{k+1})^{3/4}} \ge 2\left(\lambda_2 - \theta^{-3/4}\frac{2^{3/4}}{3}\sigma^{-1/2}\right) - \frac{2^{3/4}}{3}\sigma^{-1/2} \quad \text{a.s.}$$

Consequently,

$$\limsup_{n\to\infty} \frac{H_n}{(n\log\log n)^{3/4}} \ge 2\left(\lambda_2 - \theta^{-3/4}\frac{2^{3/4}}{3}\sigma^{-1/2}\right) - \frac{2^{3/4}}{3}\sigma^{-1/2} \quad \text{a.s.}$$

Letting $\theta \to \infty$ and $\lambda_2 \to 3^{-1}2^{3/4}\sigma^{-1/2}$ on the right-hand side gives (4.6). □

## 5. Self-intersection in high dimension

From Theorem 1.1, we have seen that the multi-dimensional case ($d \ge 2$) is different from the case $d = 1$. Here is the reason: contrary to (2.18), a concentration phenomenon appearing as

$$Q_n/\mathbb{E}Q_n \xrightarrow{p} 1$$

takes over when $d \ge 2$. The concentration also plays a role in the moderate deviations (Theorem 1.2) as $d \ge 2$. In our treatment given in Sections 2 and 3, $Q_n$ is replaced by $\mathbb{E}Q_n$ when $d \ge 2$. To justify such action, we need to show that $Q_n$ and $\mathbb{E}Q_n$ are asymptotically close enough. More precisely, our concern in this section is the central limit theorem and the exponential integrability for the renormalized self-intersection local time $Q_n - \mathbb{E}Q_n$. The case $d = 2$ has been investigated. In [20], it was proved that

$$\frac{1}{n}(Q_n - \mathbb{E}Q_n) \xrightarrow{d} (\det \Gamma)^{-1/2}\gamma_1, \tag{5.1}$$

where $\gamma_t$ is the renormalized self-intersection local times

$$\gamma_t = \iint_{0\le r<s\le t} \delta_0(W(r) - W(s))\,\mathrm{d}r\,\mathrm{d}s - \mathbb{E}\iint_{0\le r<s\le t} \delta_0(W(r) - W(s))\,\mathrm{d}r\,\mathrm{d}s, \quad t \ge 0,$$

run by a planar Brownian motion $W(t)$. In [3], it was proved that

$$\mathbb{E}\exp\left\{\frac{\lambda}{n}|Q_n - \mathbb{E}Q_n|\right\} < \infty$$

for some $\lambda > 0$. In the following discussion, we focus our attention on the case $d \geq 3$. Apart from its role in the charged polymers, the study of self-intersection local time is an important subject for its own sake. Our involvement on the integrability problems is also motivated by the recent interest [1, 2] in the large deviations for $Q_n$ in the case $d \geq 3$.

In high-dimensional cases defined by $d \geq 3$, a related object is the range $\#\{S[1,n]\}$ given by

$$\#\{S[1,n]\} = \#\{S_1, \dots, S_n\}.$$

It has been known [13, 14, 17] that

$$\frac{1}{\sqrt{n \log n}}(\#\{S[1,n]\} - \mathbb{E}\#\{S[1,n]\}) \xrightarrow{d} c_1 U \tag{5.2}$$

as $d = 3$; and

$$\frac{1}{\sqrt{n}}(\#\{S[1,n]\} - \mathbb{E}\#\{S[1,n]\}) \xrightarrow{d} c_2 U \tag{5.3}$$

as $d = 4$, where $U \sim N(0,1)$. It is now widely believed that as $d \geq 3$, $Q_n$ and $\#\{S[1,n]\}$ have very similar behaviors. In particular, we have

**Theorem 5.1.** *As $d = 3$,*

$$\frac{1}{\sqrt{n \log n}}(Q_n - \mathbb{E}Q_n) \xrightarrow{d} \lambda_1 U; \tag{5.4}$$

*as $d \geq 4$,*

$$\frac{1}{\sqrt{n}}(Q_n - \mathbb{E}Q_n) \xrightarrow{d} \lambda_2 U, \tag{5.5}$$

*where*

$$\lambda_1 = \frac{1}{\sqrt{2\pi^2 \det(\Gamma)}},$$

$$\lambda_2 = \sqrt{3G^2(0) + G(0) + 2\sum_{x \in \mathbb{Z}^d} G^3(x)},$$

$$G(x) = \sum_{k=1}^{\infty} \mathbb{P}\{S_k = x\}, \quad x \in \mathbb{Z}^d.$$

**Proof.** The proof is inspired by some ideas used in [13, 17] in the setting of the ranges. Due to similarity we only consider the case $d = 3$. Let $\{\gamma_n\}$ be a positive sequence such that

$$\gamma_n \to \infty \quad \text{and} \quad \gamma_n = \mathrm{o}(\sqrt{\log n}), \quad n \to \infty.$$

Let $0 = n_0 < n_1 < \cdots < n_{\gamma_n} = n$ be an integer partition of $[0,n]$ such that for each $1 \leq i \leq \gamma_n$, $n - i - n_{i-1} = [n\gamma_n^{-1}]$ or $[n\gamma_n^{-1}] + 1$. Then

$$Q_n = \sum_{i=1}^{\gamma_n} \sum_{n_{i-1} < j < k \leq n_i} 1_{\{S_j = S_k\}} + \sum_{i=1}^{\gamma_n - 1} \sum_{j=n_{i-1}+1}^{n_i} \sum_{k=n_i+1}^{n} 1_{\{S_j = S_k\}}. \tag{5.6}$$

For each $1 \le i \le \gamma_n - 1$,

$$\sum_{j=n_{i-1}+1}^{n_i} \sum_{k=n_i+1}^{n} 1_{\{S_j=S_k\}} \stackrel{d}{=} \sum_{j=1}^{n_i-n_{i-1}} \sum_{k=1}^{n-n_i} 1_{\{S_j=S'_k\}} \le \sum_{j=1}^{n}\sum_{k=1}^{n} 1_{\{S_j=S'_k\}},$$

where $\{S'_k\}$ is an independent copy of $\{S_k\}$. Thus,

$$\mathbb{E}\left(\sum_{i=1}^{\gamma_n-1} \sum_{j=n_{i-1}+1}^{n_i} \sum_{k=n_i+1}^{n} 1_{\{S_j=S_k\}}\right)^2 \le \gamma_n^2 \mathbb{E}\left(\sum_{j=1}^{n}\sum_{k=1}^{n} 1_{\{S_j=S'_k\}}\right)^2 = \mathrm{o}(n\log n). \tag{5.7}$$

In addition, notice that the random variables

$$\sum_{n_{i-1}<j<k\le n_i} 1_{\{S_j=S_k\}}, \quad i=1,2,\ldots,\gamma_n,$$

are independent with

$$\sum_{n_{i-1}<j<k\le n_i} 1_{\{S_j=S_k\}} \stackrel{d}{=} Q_{n_i-n_{i-1}}, \quad i=1,2,\ldots,\gamma_n.$$

By Lemma 5.1 and by (5.7),

$$\mathrm{Var}\left(\sum_{i=1}^{\gamma_n} \sum_{n_{i-1}<j<k\le n_i} 1_{\{S_j=S_k\}}\right) \sim \lambda_1^2 n\log n.$$

By Theorem 5.2, we can check the Lederberg condition. Hence,

$$\sum_{i=1}^{\gamma_n} \sum_{n_{i-1}<j<k\le n_i} 1_{\{S_j=S_k\}} \Big/ n\log n \xrightarrow{d} \lambda_1 U. \tag{5.8}$$

Finally, (5.4) follows from (5.6), (5.7) and (5.8). □

**Lemma 5.1.** *Let $\lambda_1$ and $\lambda_2$ be given in Theorem 5.1,*
*As $d=3$,*

$$\mathrm{Var}(Q_n) \sim \lambda_1^2 n\log n, \quad n\to\infty. \tag{5.9}$$

*As $d\ge 4$,*

$$\mathrm{Var}(Q_n) \sim \lambda_2^2 n, \quad n\to\infty. \tag{5.10}$$

**Proof.** Notice that

$$\begin{aligned} Q_n &= \sum_{j=1}^{n-1}\sum_{k=j+1}^{\infty} 1_{\{S_k=S_j\}} - \sum_{j=1}^{n-1}\sum_{k=n+1}^{\infty} 1_{\{S_k=S_j\}} \\ &= \sum_{j=1}^{n-1} Z_j - \sum_{j=1}^{n-1} W_j^n \quad \text{(say)}. \end{aligned} \tag{5.11}$$

Write $p_k(x) = \mathbb{P}\{S_k = x\}$ and recall that $\{S'_k\}$ is an independent copy of $\{S_k\}$. For any $i \le j \le k \le n$,

$$\begin{aligned}
\mathbb{E}W_i^n W_j^n &= \sum_{k,l=n+1}^{\infty} \mathbb{P}\{S_l - S_j = 0, S_k - S_j = S_i - S_j\} \\
&= \sum_{k,l=n-j+1}^{\infty} \mathbb{P}\{S_l = 0, S_k = S'_{j-i}\} \\
&= \sum_{x \in \mathbb{Z}^d} p_{j-i}(x) \sum_{k,l=n-j+1}^{\infty} \mathbb{P}\{S_l = 0, S_k = x\} \\
&\le \sum_{x \in \mathbb{Z}^d} p_{j-i}(x) \sum_{n-j+1 \le k \le l < \infty} \mathbb{P}\{S_l = 0, S_k = x\} \\
&\quad + \sum_{x \in \mathbb{Z}^d} p_{j-i}(x) \sum_{n-j+1 \le l \le k < \infty} \mathbb{P}\{S_l = 0, S_k = x\}.
\end{aligned} \tag{5.12}$$

For the first term on the right-hand side,

$$\begin{aligned}
&\sum_{x \in \mathbb{Z}^d} p_{j-i}(x) \sum_{n-j+1 \le k \le l < \infty} \mathbb{P}\{S_l = 0, S_k = x\} \\
&\quad = \sum_{x \in \mathbb{Z}^d} p_{j-i}(x) \sum_{n-j+1 \le k \le l < \infty} \mathbb{P}\{S_{l-k} = x\}\mathbb{P}\{S_k = x\} \\
&\quad \le C \sum_{x \in \mathbb{Z}^d} p_{j-i}(x) \sum_{n-j+1 \le k \le l < \infty} p_{l-k}(x) \frac{1}{k^{d/2}} \\
&\quad = C \sum_{n-j+1 \le k \le l < \infty} p_{l-k+j-i}(0) \frac{1}{k^{d/2}},
\end{aligned}$$

where the second step follows from the classic fact that $\sup_{x \in \mathbb{Z}^d} p_k(x) = \mathrm{O}(k^{-d/2})$.

As for the second term, a similar estimate yields that

$$\begin{aligned}
&\sum_{x \in \mathbb{Z}^d} p_{j-i}(x) \sum_{n-j+1 \le l \le k < \infty} \mathbb{P}\{S_l = 0, S_k = -x\} \\
&\quad \le C \sum_{n-j+1 \le l \le k < \infty} p_{k-l+j-i}(0) \frac{1}{l^{d/2}}.
\end{aligned}$$

Hence,

$$\begin{aligned}
\mathbb{E}W_i^n W_j^n &\le 2C \sum_{n-j+1 \le k \le l < \infty} p_{l-k+j-i}(0) \frac{1}{l^{d/2}} \\
&= 2\left(\sum_{k=n-j+1}^{\infty} p_k(0)\right)\left(\sum_{l=j-i}^{\infty} \frac{1}{l^{d/2}}\right) = \mathrm{O}((n-j)^{1-d/2}(j-i)^{1-d/2}).
\end{aligned}$$

Therefore, as $n \to \infty$,

$$\sum_{j=1}^{n} \sum_{i=1}^{j} \mathbb{E}W_i^n W_j^n = \begin{cases} \mathrm{O}(n), & d = 3, \\ \mathrm{O}((\log n)^2), & d = 4, \\ \mathrm{O}(1), & d \ge 5. \end{cases} \tag{5.13}$$

For $1 \le i \le j \le n$,

$$
\begin{aligned}
\mathrm{Cov}(Z_i, Z_j) &= \mathrm{Cov}\left(\sum_{k=j+1}^{\infty} 1_{\{S_k=S_i\}}, \sum_{k=j+1}^{\infty} 1_{\{S_k=S_j\}}\right) \\
&= \mathrm{Cov}\left(\sum_{k=1}^{\infty} 1_{\{S_k=S'_{j-i}\}}, \sum_{k=1}^{\infty} 1_{\{S_k=0\}}\right) \\
&= \sum_{k,l=1}^{\infty} \{\mathbb{P}\{S_k = S'_{j-i} S_l = 0\} - p_{k+j-i}(0)p_l(0)\} \\
&= \sum_{x\in\mathbb{Z}^d} p_{j-i}(x) \sum_{k,l=1}^{\infty} \mathbb{P}\{S_k = -x S_l = 0\} - G(0) \sum_{k=j-i+1}^{\infty} p_k(0). 
\end{aligned} \tag{5.14}
$$

Write

$$
\begin{aligned}
&\sum_{x\in\mathbb{Z}^d} p_{j-i}(x) \sum_{k,l=1}^{\infty} \mathbb{P}\{S_k = x S_l = 0\} \\
&\quad = \sum_{x\in\mathbb{Z}^d} p_{j-i}(x) \sum_{1\le k\le l<\infty} \mathbb{P}\{S_k = x S_l = 0\} \\
&\qquad + \sum_{x\in\mathbb{Z}^d} p_{j-i}(x) \sum_{1\le l<k<\infty} \mathbb{P}\{S_k = x S_l = 0\}.
\end{aligned}
$$

We have

$$
\begin{aligned}
&\sum_{x\in\mathbb{Z}^d} p_{j-i}(x) \sum_{1\le k\le l<\infty} \mathbb{P}\{S_k = x S_l = 0\} \\
&\quad = \sum_{x\in\mathbb{Z}^d} p_{j-i}(x) \sum_{1\le k\le l<\infty} p_k(x)p_{l-k}(x) \\
&\quad = \sum_{x\in\mathbb{Z}^d} p_{j-i}(x)G(x) \sum_{l=0}^{\infty} p_l(x) \\
&\quad = p_{j-i}(0)G(0)(1+G(0)) + \sum_{x\ne 0} p_{j-i}(x)G^2(x)
\end{aligned}
$$

and

$$
\begin{aligned}
&\sum_{x\in\mathbb{Z}^d} p_{j-i}(x) \sum_{1\le l<k<\infty} \mathbb{P}\{S_k = x S_l = 0\} \\
&\quad = \sum_{x\in\mathbb{Z}^d} p_{j-i}(x) \sum_{1\le l<k<\infty} p_{k-l}(x)p_l(0) \\
&\quad = \sum_{1\le l<k<\infty} p_{k-l+j-i}(0)p_l(0) = G(0) \sum_{k=j-i+1}^{\infty} p_k(0).
\end{aligned}
$$

In summary of the argument since (5.14),

$$
\mathrm{Cov}(Z_i, Z_j) = p_{j-i}(0)G(0)(1+G(0)) + \sum_{x\ne 0} p_{j-i}(x)G^2(x).
$$

Consequently,

$$\begin{aligned}\mathrm{Var}\left(\sum_{i=1}^n Z_i\right) &= \sum_{i=1}^n \mathrm{Var}(Z_i) + 2\sum_{1\le i<j\le n}\mathrm{Cov}(Z_i,Z_j)\\ &= nG(0)(1+G(0)) + 2G(0)(1+G(0))\sum_{j=1}^{n-1}G_{n-j}(0) + 2\sum_{x\neq 0}G^2(x)\sum_{i=1}^{n-1}G_{n-i}(x)\\ &= nG(0)(1+G(0)) + 2G(0)\sum_{j=1}^{n-1}G_j(0) + 2\sum_{x\in\mathbb{Z}^d}G^2(x)\sum_{j=1}^{n-1}G_j(x),\end{aligned}\tag{5.15}$$

where

$$G_j(x) = \sum_{i=1}^j p_i(x).$$

When $d\ge 4$,

$$\sum_{x\in\mathbb{Z}^d} G^3(x) < \infty.$$

By (5.15)

$$\mathrm{Var}\left(\sum_{i=1}^n Z_i\right) \sim n\left\{G(0)(1+G(0)) + 2G^2(0) + 2\sum_{x\in\mathbb{Z}^d}G^3(x)\right\},\quad n\to\infty.$$

By (5.11) and (5.13), this implies (5.10).

We now consider the case $d=3$. We use the fact that (p. 308, [21])

$$G(x) \sim (2\pi)^{-1}\det(\Gamma)^{-1/2}\langle x,\Gamma^{-1}x\rangle^{-1/2},\quad |x|\to\infty.\tag{5.16}$$

By (5.16),

$$\sum_{\langle x,\Gamma^{-1}x\rangle > j} G(x)^2 G_j(x) = \mathrm{O}\left(\frac{1}{j}\sum_{x\in\mathbb{Z}^d}G_j(x)\right) = \mathrm{O}(1),\quad j\to\infty.\tag{5.17}$$

In addition,

$$\begin{aligned}&\sum_{\langle x,\Gamma^{-1}x\rangle\le j} G(x)^2(G(x)-G_j(x))\\ &\quad\le (G(0)-G_j(0))\sum_{\langle x,\Gamma^{-1}x\rangle\le j}G^2(x)\\ &\quad\le C(G(0)-G_j(0))\sum_{\langle x,\Gamma^{-1}x\rangle\le j}\frac{1}{1+|x|^2} = \mathrm{O}(1),\quad j\to\infty,\end{aligned}\tag{5.18}$$

and

$$\begin{aligned}\sum_{\langle x,\Gamma^{-1}x\rangle\le j} G^3(x) &\sim \frac{1}{(2\pi)^3\det(\Gamma)^3}\int_{\{1\le\langle x,\Gamma^{-1}x\rangle\le j\}}\langle x,\Gamma^{-1}x\rangle^{-3/2}\,\mathrm{d}x\\ &= \frac{1}{(2\pi)^3\det(\Gamma)}\int_{\{1\le|y|\le\sqrt{j}\}}|y|^{-3}\,\mathrm{d}y = \frac{1}{(2\pi)^2\det(\Gamma)}\log j,\quad j\to\infty.\end{aligned}\tag{5.19}$$

Combining (5.17), (5.18) and (5.19) gives

$$\sum_{x\in\mathbb{Z}^d} G^2(x)G_j(x) \sim \frac{1}{(2\pi)^2\det(\Gamma)}\log j, \quad j\to\infty.$$

By (5.14),

$$\mathrm{Var}\left(\sum_{i=1}^n Z_i\right) \sim \frac{2}{(2\pi)^2\det(\Gamma)} n\log n, \quad n\to\infty,$$

which, together with (5.11) and (5.13), implies (5.9). □

We now investigate the integrability of $Q_n$. Write

$$J_n = \sum_{j=1}^n\sum_{k=1}^n 1_{\{S_j=S'_k\}},$$

where $\{S'_n\}$ is an independent copy of $\{S'_n\}$. $J_n$ is known as the intersection local time between two independent trajectories.

**Lemma 5.2.** *As $d\ge 3$, there is a constant $C_d>0$ such that*

$$\mathbb{E}J_n^m \le C_d^m (m!)^{3/2} n^{m/2}, \quad m,n=1,2,\ldots. \tag{5.20}$$

**Proof.** Recall the fact (p. 3282, [6])

$$\mathbb{E}J_n^m \le (m!)(1+\mathbb{E}J_n)^m, \quad m,n=1,2,\ldots, \tag{5.21}$$

and the fact that

$$\mathbb{E}J_n = \begin{cases} \mathrm{O}(\sqrt{n}), & d=3,\\ \mathrm{O}(\log n), & d=4,\\ \mathrm{O}(1), & d\ge 5.\end{cases} \tag{5.22}$$

A trivial and rough summarization of (5.22) gives that $\mathbb{E}J_n = \mathrm{O}(\sqrt{n})$ as $d\ge 3$. By (5.21), we obtain a weaker version of (5.20):

$$\mathbb{E}J_n^m \le C^m (m!)^2 n^{m/2}, \quad m,n=1,2,\ldots. \tag{5.23}$$

To strengthen (5.23) into (5.20), recall (Theorem 5.1, of [6]) that for any integers $m,n_1,\ldots,n_a\ge 1$,

$$(\mathbb{E}J_n^m)^{1/2} \le \sum_{\substack{k_1+\cdots+k_a=m\\ k_1,\ldots,k_a\ge 0}} \frac{m!}{k_1!\cdots k_a!}(\mathbb{E}J_{n_1}^{k_1})^{1/2}\cdots(\mathbb{E}J_{n_a}^{k_m})^{1/2}. \tag{5.24}$$

We first consider the case $n\ge m$. Write $l(m,n)=1+[\frac{n}{m}]$. By (5.23) there is a $C>0$ independent of $m$ and $n$ such that

$$\mathbb{E}J_{l(m,n)}^m \le C^m (m!)^2\left(\frac{n}{m}\right)^{m/2}.$$

Taking $a = m$ in (5.24) gives

$$(\mathbb{E}J_n^m)^{1/2} \le \sum_{\substack{k_1+\cdots+k_m=m\\ k_1,\ldots,k_m\ge 0}} \frac{m!}{k_1!\cdots k_m!}(\mathbb{E}J_{l(m,n)}^{k_1})^{1/2}\cdots(\mathbb{E}J_{l(m,n)}^{k_m})^{1/2}$$
$$\le m!C^{m/2}\left(\frac{n}{m}\right)^{m/4} \sum_{\substack{k_1+\cdots+k_a=m\\ k_1,\ldots,k_m\ge 0}} 1 = m!C^{m/2}m^{-m/4}n^{m/4}\binom{2m-1}{m}.$$

Since $m^m \ge m!$ and

$$\binom{2m-1}{m} \le 4^m$$

we have established (5.20) with $C_d = 8C$ in the case $n \ge m$.

As for the case $n < m$, the trivial fact $J_n \le n^2$ leads to the following trivial bound,

$$\mathbb{E}J_n^m \le n^{2m} \le m^{3m/2}n^{m/2} \le C^m(m!)^{3/2}n^{m/2},$$

where the last step follows from the Stirling formula. □

As for the exponential integrability of the renormalized self-intersection local time $Q_n - \mathbb{E}Q_n$, we have the following theorem.

**Theorem 5.2.** *As* $d = 3$,

$$\sup_n \mathbb{E}\exp\left\{\frac{\theta}{\sqrt[3]{n\log n}}|Q_n - \mathbb{E}Q_n|^{2/3}\right\} < \infty \quad \textit{for every } \theta > 0. \tag{5.25}$$

*As* $d \ge 4$,

$$\sup_n \mathbb{E}\exp\left\{\frac{\theta}{\sqrt[3]{n}}|Q_n - \mathbb{E}Q_n|^{2/3}\right\} < \infty \quad \textit{for some } \theta > 0. \tag{5.26}$$

**Proof.** The proof given here is radically different from the approach used in the case $d = 2$ (Lemma 2.3 in [3]) where the treatment is Le Gall–Varadhan's triangular approximation. Due to similarity, we only consider $d = 3$. We first prove that for any integer $m \ge 1$,

$$\mathbb{E}|Q_n - \mathbb{E}Q_n|^m = O((n\log n)^{m/2}), \quad n \to \infty. \tag{5.27}$$

We carry out induction on $m$. By Lemma 5.2, (5.27) holds as $m = 1, 2$. We let $m \ge 3$, assume that it is true for all $1 \le j \le m - 1$ and prove it is true for $m$.

Given $n$, write $n_1 = [n/2]$, $n_2 = n - n_1$,

$$Q'_{n_2} = \sum_{n_1+1\le j<k\le n} 1_{\{S_j=S_k\}} \quad \text{and} \quad \bar{J}_n = \sum_{j=1}^{n_1}\sum_{k=n_1+1}^{n} 1_{\{S_j=S_k\}}.$$

By Lemma 5.2 there is $C_0 > 0$ such that

$$\mathbb{E}|\bar{J}_n - \mathbb{E}\bar{J}_n|^m \le C_0^m(m!)^{3/2}n^{m/2}m, \quad n = 1, 2, \ldots. \tag{5.28}$$

By independence between $Q_{n_1}$ and $Q'_{n_2}$,

$$\begin{aligned}
&(\mathbb{E}|Q_n - \mathbb{E}Q_n|^m)^{1/m} \\
&\quad \le \{\mathbb{E}(|Q_{n_1} - \mathbb{E}Q_{n_1}| + |Q'_{n_2} - \mathbb{E}Q'_{n_2}|)^m\}^{1/m} + \{\mathbb{E}|\bar{J}_n - \mathbb{E}\bar{J}_n|^m\}^{1/m} \\
&\quad \le \left\{\sum_{j=0}^m \binom{m}{j} \mathbb{E}|Q_{n_1} - \mathbb{E}Q_{n_1}|^j \mathbb{E}|Q'_{n_2} - \mathbb{E}Q'_{n_2}|^{m-j}\right\}^{1/m} \\
&\qquad + C_0\sqrt{n}(m!)^{3/(2m)}. \qquad (5.29)
\end{aligned}$$

Combining the induction assumption and (5.29),

$$\begin{aligned}
&(\mathbb{E}|Q_n - \mathbb{E}Q_n|^m)^{1/m} \\
&\quad \le \{\mathrm{O}((n\log n)^{m/2}) + \mathbb{E}|Q_{n_1} - \mathbb{E}Q_{n_1}|^m + \mathbb{E}|Q'_{n_2} - \mathbb{E}Q'_{n_2}|^m\}^{1/m} + \mathrm{O}(\sqrt{n}).
\end{aligned}$$

Write

$$\alpha_k = \sup\{(2^k \log 2^k)^{-1/2}(\mathbb{E}|Q_n - \mathbb{E}Q_n|^m)^{1/m}; 2^k \le n \le 2^{k+1}\}.$$

Then,

$$\alpha_{k+1} \le \{\mathrm{O}(1) + 2^{-(m-2)/2}\alpha_k^m\}^{1/m} + \mathrm{o}(1) \le 2^{-(m-2)/(2m)}\alpha_k + \mathrm{O}(1), \quad k \to \infty.$$

By the fact that $m \ge 3$ one can see that the sequence $\alpha_k$ is bounded. We have proved (5.27).

We now claim that there is $C > 0$ such that

$$\mathbb{E}|Q_n - \mathbb{E}Q_n|^m \le C^m(m!)^{3/2}(n\log n)^{m/2}, \quad m, n = 1, 2, \ldots. \qquad (5.30)$$

Indeed, take $m_0$ sufficiently large so that

$$(1 - 2^{-(m-2)/(2m)})^{-1}\frac{(m-1)^{1/m}}{\sqrt{2}} \le \frac{1}{2}, \qquad 1 - 2^{-(m-2)/(2m)} \ge \frac{1}{4}$$

for all $m \ge m_0$. By (5.27), there is a constant $C > 0$ such that for all $j = 1, \ldots, m_0$,

$$\mathbb{E}|Q_n - \mathbb{E}Q_n|^j \le C^j(j!)^{3/2}(n\log n)^{j/2}, \quad n = 1, 2, \ldots.$$

We may assume that $C \ge 8C_0$. (Recall that $C_0$ is given in (5.28).) By induction (on $m$), all we have to prove is that for any $m \ge m_0$, if

$$\mathbb{E}|Q_n - \mathbb{E}Q_n|^j \le C^j(j!)^{3/2}(n\log n)^{j/2}, \quad n = 1, 2, \ldots, \qquad (5.31)$$

for every $j = 1, \ldots, m-1$ then

$$\mathbb{E}|Q_n - \mathbb{E}Q_n|^m \le C^m(m!)^{3/2}(n\log n)^{m/2}, \quad n = 1, 2, \ldots. \qquad (5.32)$$

From (5.29) and (5.31) we have

$$\begin{aligned}
(\mathbb{E}|Q_n - \mathbb{E}Q_n|^m)^{1/m} \le \Bigg\{&2^{-m/2}C^m(n\log n)^{m/2}\sum_{j=1}^{m-1}\binom{m}{j}(j!)^{3/2}((m-j)!)^{3/2} \\
&+ \mathbb{E}|Q_{n_1} - \mathbb{E}Q_{n_1}|^m + \mathbb{E}|Q'_{n_2} - \mathbb{E}Q'_{n_2}|^m\Bigg\}^{1/m} + C_0\sqrt{n}(m!)^{3/(2m)}.
\end{aligned}$$

Notice that

$$\sum_{j=1}^{m-1}\binom{m}{j}(j!)^{3/2}((m-j)!)^{3/2}=m!\sum_{j=1}^{m-1}(j!)^{1/2}((m-j)!)^{1/2}\le(m!)^{3/2}(m-1)$$

and that by (5.27)

$$\beta_m\equiv\sup_{n\ge1}\{(n\log n)^{-m/2}\mathbb{E}|Q_n-\mathbb{E}Q_n|^m\}<\infty.$$

We have

$$\begin{aligned}\beta_m^{1/m}&\le\left\{\frac{m-1}{2^{m/2}}(m!)^{3/2}C^m+2^{-(m-2)/2}\beta_m\right\}^{1/m}+C_0(m!)^{3/(2m)}\\&\le C\frac{(m-1)^{1/m}}{\sqrt{2}}(m!)^{3/(2m)}+2^{-(m-2)/(2m)}\beta_m^{1/m}+C_0(m!)^{3/(2m)}.\end{aligned}$$

Hence,

$$\begin{aligned}\beta_m^{1/m}&\le(1-2^{-(m-2)/(2m)})^{-1}\left\{\frac{(m-1)^{1/m}}{\sqrt{2}}C+C_0\right\}(m!)^{3/(2m)}\\&\le\left(\frac{1}{2}C+4C_0\right)(m!)^{3/(2m)}\le C(m!)^{3/(2m)}.\end{aligned}$$

Hence, (5.27) holds. By (5.27) and the Taylor expansion there is $\theta_0>0$ such that

$$\sup_n\mathbb{E}\exp\left\{\frac{\theta_0}{\sqrt[3]{n\log n}}|Q_n-\mathbb{E}Q_n|^{2/3}\right\}<\infty. \tag{5.33}$$

It remains to extend (5.33) to any $\theta>0$. Indeed, for any $\theta>\theta_0$, one can find an integer $l$ such that for any $n$ there is an integer partition $0=n_0<n_1,\dots,<n_l=n$ such that $n_i-n_{i-1}<n(\theta_0/\theta)^3$ $(i=1,\dots,l)$. Write

$$Q_n=\sum_{i=1}^{l}\sum_{n_{i-1}<j,k\le n_i}1_{\{S_j=S_k\}}+\sum_{i=2}^{l}\sum_{j=1}^{n_{i-1}}\sum_{k=n_{i-1}+1}^{n_i}1_{\{S_j=S_k\}}. \tag{5.34}$$

Notice that for each $2\le i\le l$,

$$\sum_{j=1}^{n_{i-1}}\sum_{k=n_{i-1}+1}^{n_i}1_{\{S_j=S_k\}}\overset{d}{=}\sum_{j=1}^{n_{i-1}}\sum_{k=1}^{n_i-n_{i-1}}1_{\{S_j=S'_k\}}\le J_n.$$

By Lemma 5.2 and the Taylor expansion there is $\lambda>0$ such that

$$\sup_{n\ge1}\mathbb{E}\exp\left\{\frac{\lambda}{\sqrt[3]{n}}\left(\sum_{j=1}^{n_{i-1}}\sum_{k=n_{i-1}+1}^{n_i}1_{\{S_j=S_k\}}\right)^{2/3}\right\}<\infty.$$

By the triangular inequality and the Hölder inequality, therefore,

$$\sup_{n\ge1}\mathbb{E}\exp\left\{\frac{\theta}{\sqrt[3]{n\log n}}\left(\sum_{i=2}^{l}\sum_{j=1}^{n_{i-1}}\sum_{k=n_{i-1}+1}^{n_i}1_{\{S_j=S_k\}}\right)^{2/3}\right\}<\infty.$$

Consequently, by the Jensen inequality

$$\sup_{n\geq 1}\mathbb{E}\exp\left\{\frac{\theta}{\sqrt[3]{n\log n}}\left(\mathbb{E}\sum_{i=2}^{l}\sum_{j=1}^{n_{i-1}}\sum_{k=n_{i-1}+1}^{n_i}1_{\{S_j=S_k\}}\right)^{2/3}\right\}<\infty.$$

By independence and by the triangular inequality

$$\begin{aligned}&\mathbb{E}\exp\left\{\frac{\theta}{\sqrt[3]{n\log n}}\left(\sum_{i=1}^{l}\left|\sum_{n_{i-1}<j,k\leq n_i}1_{\{S_j=S_k\}}-\mathbb{E}\sum_{n_{i-1}<j,k\leq n_i}1_{\{S_j=S_k\}}\right|\right)^{2/3}\right\}\\ &\quad\leq\prod_{i=1}^{l}\mathbb{E}\exp\left\{\frac{\theta}{\sqrt[3]{n\log n}}|Q_{n_i-n_{i-1}}-\mathbb{E}Q_{n_i-n_{i-1}}|^{2/3}\right\}\\ &\quad\leq\prod_{i=1}^{l}\mathbb{E}\exp\left\{\frac{\theta_0}{\sqrt[3]{(n_i-n_{i-1})\log(n_i-n_{i-1})}}|Q_{n_i-n_{i-1}}-\mathbb{E}Q_{n_i-n_{i-1}}|^{2/3}\right\}.\end{aligned}$$

Hence, the desired (5.25) follows from (5.33) and (5.34). □

By slightly modifying our argument, we can prove that Theorem 5.2 holds also for the range. We include this in our paper for future reference.

**Theorem 5.3.** *As* $d=3$,

$$\sup_{n}\mathbb{E}\exp\left\{\frac{\theta}{\sqrt[3]{n\log n}}|\#\{S[1,n]\}-\mathbb{E}\#\{S[1,n]\}|^{2/3}\right\}<\infty\quad\text{for every }\theta>0. \tag{5.35}$$

*As* $d\geq 4$,

$$\sup_{n}\mathbb{E}\exp\left\{\frac{\theta}{\sqrt[3]{n}}|\#\{S[1,n]\}-\mathbb{E}\#\{S[1,n]\}|^{2/3}\right\}<\infty\quad\text{for some }\theta>0. \tag{5.36}$$

**Proof.** The argument used here is essentially the same as the one for $Q_n$. To carry it through, (5.21) needs to be replaced by Lemma 3.1 in [17], and (5.24) needs to be replaced by Lemma 6 in [7]. The rest of the proof follows an obvious modification of the argument for Theorem 5.2. □

## Acknowledgments

The author thanks anonymous referee(s) for constructive comments and suggestions.